\documentclass{amsart}

\usepackage{xspace,amssymb,amsfonts,euscript,eufrak,mathrsfs}

\usepackage{amsthm,amsmath}
\usepackage{url}
\usepackage{graphicx}
 \usepackage{epstopdf}
\usepackage{xcolor}
\usepackage{color}
\usepackage[all]{xy}
\usepackage{geometry}

\usepackage{tikz}
\usepackage{pgflibraryshapes}
\usetikzlibrary{arrows}
\usetikzlibrary{shapes}
\usetikzlibrary{decorations.pathmorphing,patterns,calc}
\usetikzlibrary{automata}
\tikzstyle{model_node}=[circle,draw=black,fill=tugreen!50,minimum size=50pt]
\tikzstyle{tool}=[draw,double,rounded corners,inner sep=10pt]
\tikzstyle{aut_node}=[circle,draw=black,fill=black,inner sep=0pt,minimum size=7pt,font=\footnotesize]

\renewcommand{\subsectionmark}[1]{\markright{}{}}

\newcommand{\g}{\mathfrak{g}}

\newcommand{\al}{\alpha}
\newcommand{\als}{\alpha_s}
\newcommand{\alt}{\alpha_t}
\newcommand{\BMP}{\mathscr{B}}

\newcommand{\E}{\mathcal{E}}
\newcommand{\F}{\mathscr{F}}
\newcommand{\He}{\mathbf{H}}
\renewcommand{\H}{\mathcal{H}}
\newcommand{\HJ}{\mathcal{H}^J}

\newcommand{\sHJ}{{}^s\mathcal{H}^J}
\newcommand{\I}{\mathcal{I}}

\newcommand{\MeJ}{\mathbf{M}^J}

\newcommand{\MG}{\mathcal{G}}

\newcommand{\MGJ}{\mathcal{G}^J}

\newcommand{\SR}{\mathcal{S}}

\newcommand{\ShMGk}{\textbf{Sh}(\MG)}
\newcommand{\Sh}[1]{\textbf{Sh}({#1})}
\newcommand{\ShZ}{\mathscr{Z}}
\newcommand{\T}{\mathcal{T}}

\newcommand{\tle}{\trianglelefteq}

\newcommand{\V}{\mathcal{V}}
\newcommand{\W}{\mathcal{W}}

\newcommand{\Z}{\mathcal{Z}}
\newcommand{\ZJ}{\mathcal{Z}^J}

\newcommand{\sZJ}{{}^s\mathcal{Z}^J}
\newcommand{\msZJ}{{}^{-s}\mathcal{Z}^J}

\newcommand{\sth}{{}^s\theta}
\newcommand{\sthJ}{{}^s\theta^J}
\newcommand{\ssig}{{}_s\sigma}
\newcommand{\0}{\emptyset}

\renewcommand{\l}{\ell}

\newcommand{\Hom}{\text{Hom}}
\newcommand{\im}{\text{im}}
\newcommand{\rk}{\underline{\text{rk}}\,}

\newtheorem{theor}{Theorem}[section]
\newtheorem{prop}[theor]{Proposition}
\newtheorem*{theor*}{Theorem}
\newtheorem{cor}[theor]{Corollary}
\newtheorem{lem}[theor]{Lemma}

\theoremstyle{definition}
\newtheorem{rem}[theor]{Remark}

\newtheorem{expl}[theor]{Example}
\newtheorem{defin}[theor]{Definition}

\begin{document}
\title{Categorification of a parabolic Hecke module via sheaves on moment graphs}
\author{Martina Lanini}
\address{
Department of Mathematics and Statistics, University of Melbourne, Parkville VIC 3010 Australia}
\email{martina.lanini@unimelb.edu.au}
\maketitle
\vspace{-9mm}
\begin{abstract}
We investigate certain categories, associated by Fiebig with the geometric representation of a Coxeter system, via sheaves on Bruhat graphs. 
We modify Fiebig's definition of translation functors in order to extend it to the singular setting and use it to categorify a parabolic Hecke module. As an application we obtain a combinatorial description of indecomposable projective objects of (truncated) non-critical singular blocks of (a deformed version of) category $\mathcal{O}$, using indecomposable special modules over the structure algebra of the corresponding Bruhat graph.
\end{abstract}

\section{Introduction}

A typical problem in the representation theory of Kac-Moody algebras is to understand the composition series of standard objects in the corresponding category $\mathcal{O}$ of Bernstein, Gelfand and Gelfand (cf. \cite{BGG}).
 In the case of a standard object lying in a regular block, this question is the core of the Kazhdan-Lusztig theory, and the answer is known to be given by the Kazhdan-Lusztig polynomials evaluated at the identity.
 If we consider a singular block, we only have to replace these polynomials by their parabolic analogue. 
In the case of a principal block, this fact was conjectured by Kazhdan and Lusztig in \cite{KL} and proved in several steps in \cite{KL80}, \cite{BeBe}, \cite{BK}. A fundamental role in the proof of the Kazhdan-Lusztig
 conjecture was played by the geometric interpretation of the problem in terms of perverse sheaves and intersection cohomology complexes. In particular, one could study certain properties of the Hecke algebra in the category
 of equivariant  perverse sheaves on the corresponding flag variety. 

An alternative way to attack the Kazhdan-Lusztig conjecture is via Soergel bimodules, which provide a combinatorial realisation of projective objects in category $\mathcal{O}$. The combinatorial description of indecomposable projective objects we present in this paper is an analogue of  Soergel's combinatorial contruction (introduced at first  for finite dimensional Lie algebras in \cite{Soe90}). The Soergel bimodule approach to the Kazhdan-Lusztig conjecture recently led to an algebraic proof of it  by Elias and Williamson (cf. \cite{EW}).

The procedure of considering a complicated object, such as a category, in order to understand a simpler one is motivated by the fact that 
the extra structure can provide us with new tools and allow us to prove and hopefully  generalise certain phenomena that are difficult to address directly.

In \cite{Deo87}, Deodhar associated with any Coxeter system $(\W, \SR)$ and any subset of the set of simple reflections $J\subseteq \SR$ the parabolic Hecke module $\MeJ$.  The aim of this paper is to give a categorification of this module, for any $J$  generating a 
finite subgroup. 

 We have followed the definition of categorification of $\MeJ$ as in \cite[Remark 7.8]{MS1}, which is actually a weak categorification. This could be strengthened to a proper categorification by presenting the result as a 2-representation of some 2-category (see \cite[Sections 1-3]{M} for various levels of categorification  and Remark \ref{rem_2cat} of this paper for a more precise statement). In \cite{MS2}, the authors properly categorify induced cell modules (in the finite case), which is a huge step outside the parabolic Hecke module (the latter being just a special case).

If $\W$ is a Weyl group, there is a partial flag variety $Y$ corresponding to the set $J$, equipped with an action of 
a maximal torus $T$, and, as for the regular case, one possible categorification is given by the category of $B$-equivariant perverse 
sheaves on $Y$. Our goal is to describe a general categorification, which can be defined also in the case in which there is no geometry 
available. In order to do this, our main tools will be  Bruhat moment graphs and sheaves on them. We will see how these objects come 
naturally into the picture.

Moment graphs appeared for the first time in \cite{GKM} as 1-skeletons of  actions of  tori  on  complex algebraic 
varieties. In particular, Goresky, Kottwitz and MacPherson were able to  describe explicitly the equivariant cohomology of these
 varieties  using only the data encoded in the underlying moment graphs. 
Inspired by this result, Braden and MacPherson (cf. \cite{BM01})  could study the equivariant intersection cohomology of a complex
 algebraic variety equipped with a Whitney stratification, stable with respect to the torus action. In order to do so, they introduced the notion of sheaves on moment graphs and, in particular, of \emph{canonical sheaves}. We will refer to this class of sheaves as \emph{Braden-MacPherson(BMP)-sheaves}.

Even if moment graphs arose originally from geometry, Fiebig observed that it is possible to give an axiomatic definition of them 
(cf. \cite{Fie08b}). In particular, he associated a moment graph to any Coxeter datum $(\W, \SR, J)$ as above and, in the case
 of $J=\emptyset$, he used it to give an alternative construction of Soergel's category of bimodules associated to a reflection faithful
 representation of $(\W,\SR)$ (cf. \cite{Fie08b}). (We refer the reader to \cite{W} for the singular version of Soergel's bimodules.) The indecomposable objects of the category defined by Fiebig are precisely the \emph{BMP}-sheaves, that, if $\W$ is a Weyl group,  are related to the  intersection cohomology complexes, the simple objects in the category of perverse sheaves.  A fundamental step in Fiebig's realisation of this category were translation functors, whose definition we extend to  the parabolic setting (see \S\ref{sssec_TransLeft}).
 
The paper is organised as follows. 

In Section 2 we recall the definition of the parabolic Hecke module $\MeJ$ and the fact that it is the unique free $\mathbb{Z}[v,v^{-1}]$-module
 having rank $|\W/\langle J \rangle|$ equipped with a certain structure of a module over the Hecke algebra $\He$. This structure is described 
in terms of the action of the Kazhdan-Lusztig basis elements $\underline{H}_s$, for $s\in \SR$.
 Then by a \emph{categorification} of $\MeJ$ (as in \cite[Remark 7.8]{MS1}) we  mean a category $\mathscr{C}$, which is exact in the sense of Quillen (cf. \cite{Qu}), together with an autoequivalence $G$ and 
exact functors  $\{F_s\}_{s\in\SR}$, that provide the Grothendieck group $[\mathscr{C}]$  with the structure of a $\mathbb{Z}[v,v^{-1}]$-module and $\He$-module, such that there exists an isomorphism from $[\mathscr{C}]$ to the parabolic module, satisfying certain 
compatibility conditions with these functors coming from the defining properties of $\MeJ$ (see Definition \ref{ssec_categorification}).

In the third section we introduce the objects we will be dealing with in the rest of the paper. In particular, we review basic concepts of the theory of moment graphs and sheaves on them.

 Section  4 is about $\mathbb{Z}$-graded modules over $\ZJ$, the structure algebra of a parabolic Bruhat graph. In particular, for any
 $s\in \SR$, we define the translation functor $\sth$ and define the category $\HJ$ of \emph{special $\ZJ$-modules}. 
By definition, this category is stable under the shift in degree, that we denote by $\langle\cdot \rangle$, and under $\sth$ for all $s\in\SR$.

 In Section 5 we study certain subquotients
 of objects in $\HJ$ and this allows us to define an exact structure on $\HJ$ and hence to state our main theorem:

\vspace{3mm}
\textbf{Theorem \ref{Thm_catHM}} \textit{The category $\H^J$ special $\ZJ$-modules together with the shift in degree $\langle -1 \rangle$ and (shifted) translation functors is
 a categorification of the parabolic Hecke module $\MeJ$.}
\vspace{3mm}

Section 6 is devoted to the proof of this theorem. First of all, we show that $\sth\circ\langle 1 \rangle$ is an exact functor (Lemma \ref{sth_exact}). Secondly, we define the character map $h^J:[\HJ]\rightarrow \MeJ$ and  prove that the functors  $\langle -1 \rangle$ and  $\sth\circ \langle 1\rangle$, $s\in \SR$,  satisfy the desired 
compatibility condition (Proposition \ref{prop_left_mult}). We conclude then by showing that the character map is an isomorphism of $\mathbb{Z}[v,v^{-1}]$-modules (Lemma \ref{Lem_surj_hJ} and Lemma \ref{Lem_inj_hJ}).

Section 7 is about the  categorification of a certain injective map  of 
$\He$-modules $i:\MeJ\hookrightarrow \He$, which allows us to see  the category $\HJ$ as a subcategory of $\H^{\0}$. More precisely, we define an exact functor $I:\HJ\rightarrow \H^{\emptyset}$ such that
the following diagram commutes:
\begin{displaymath}
 \begin{xymatrix}{
[\HJ]\ar[d]_{h^J}\ar@{^{(}->}[r]^{[I]}&[\H^{\0}]\ar[d]^{h^{\0}}\\
\MeJ \ar@{^{(}->}[r]_{ i }&\He &}
\end{xymatrix}
\end{displaymath}

In order to construct and investigate the functor $I$, we give a realisation of $\HJ$ via \emph{BMP}-sheaves (Proposition \ref{prop_HJ_BMP}) and then use Fiebig's idea of interchanging global and local viewpoints (cf. \cite{Fie08b}).

In the last section we discuss briefly the relationship between $\HJ$ and non-critical blocks of an equivariant version of category $\mathcal{O}$ for symmetrisable Kac-Moody algebras. In particular, we show that the indecomposable projective objects of a truncated, non-critical block  $\mathcal{O}_{R, \Lambda^{\leq \nu}}$ are combinatorially described by indecomposable modules in $\HJ$, with $J$  depending on $\Lambda$ (Proposition \ref{prop_eqvtRepnTh}).

\section{Hecke modules}\label{sec_HeckeM}

Here we recall some classical constructions, following \cite{Soe97}. We close the section by defining the concept of categorification of the parabolic Hecke module $\MeJ$. 

\subsection{Hecke algebra}
The Hecke algebra associated to a Coxeter system $(\W, \SR)$ is nothing but a quantisation of the group ring $\mathbb{Z}[\W]$. Let $\leq$ be the Bruhat order on $\W$ and $\l:\W\rightarrow \mathbb{Z}$ be the length function associated to $\SR$. Denote by $\mathcal{L}:=\mathbb{Z}[v,v^{-1}]$ the ring of Laurent polynomials in the variable $v$ over $\mathbb{Z}$.

\begin{defin}The Hecke algebra $\He=\He(\W, \SR)$ is the free $\mathcal{L}$-module having basis $\{H_x\,|\, x\in \W\}$, subject to the following relations: 
\begin{equation}\label{eqn_Hecke}H_sH_x=\left\{
\begin{array}{lcr}
H_{sx} & \text{\small{if}} & sx>s, \\
(v^{-1}-v)H_{x}+H_{sx}  &  \text{\small{if}} & sx<x,
\end{array}\right.\qquad \hbox{\text{for}}\quad x\in \W, \ s\in\SR.
\end{equation}
\end{defin}

It is well known that this defines an associative $\mathcal{L}$-algebra (cf. \cite{Humph}).

It is easy to verify that $H_x$ is invertible for any $x\in \W$, and this allows us to  define an involution on $\He$; i.e. the unique ring homomorphism $\stackrel{\_\!\_\!\_}{\ \cdot\  }\  :\He\rightarrow\He$ such that $\overline{v}=v^{-1}$ and $\overline{H_x}=(H_{x^{-1}})^{-1}$.

In \cite{KL} Kazhdan and Lusztig showed the existence of another basis $\{\underline{H}_x\}$ for $\He$, the so--called \emph{Kazhdan-Lusztig  basis}, that they used to define complex representations of the Hecke algebra and hence of the Coxeter group.
 The entries of the change of basis matrix are given by a family of polynomials in $\mathbb{Z}[v]$, which are called  \emph{Kazhdan-Lusztig  polynomials}. There are many different normalisations of this basis appearing in the literature. The one we adopt, following \cite{Soe97},  is determined by Theorem \ref{Deo87} (see  Remark \ref{rem_parab_regularHM}).

\subsubsection{Parabolic Hecke modules}\label{sssec_HeckeKL}
 In \cite{Deo87} Deodhar generalised this construction to the parabolic setting in the following way. Let $\W,\SR$ and $\He$ be as above. Fix a subset $J\subseteq \SR$ and denote by $\W_J=\langle J \rangle$ the subgroup of
 $\W$ generated by $J$. Since $(\W_J, J)$ is also a Coxeter system, it makes sense to consider its Hecke algebra $\He_J=\He(\W_J, J)$. 

For any simple reflection $s\in \SR$,  the element $H_s$ satisfies the quadratic relation $(H_s)^2=(v^{-1}-v)H_{s}+H_{e}$; that is,
  $(H_s+v)(H_s-v^{-1})=0$. If $u\in \{v^{-1},-v\}$, we may define a map of $\mathcal{L}$-modules $\varphi_u:\He_J\rightarrow \mathcal{L}$
by $H_s\mapsto u$. In this way, $\mathcal{L}$ is endowed with the structure of a $\He_J$-bimodule, which we denote by $\mathcal{L}(u)$.

The parabolic Hecke modules are defined as $\mathbf{M}^J:=\He\otimes_{\He_J}\mathcal{L}(v^{-1})$ and $\mathbf{N}^{J}:=\He\otimes_{\He_J}\mathcal{L}(-v)$. 
As in the Hecke algebra case, it is possible to define an involutive automorphism of these modules. Namely,
\begin{equation}\begin{array}{rcl}\stackrel{\_\!\_\!\_}{\ \cdot\  }\ :\ \He \otimes_{\He_J}\mathcal{L}(u)&\rightarrow&\He \otimes_{\He_J}\mathcal{L}(u)\\
H\otimes a&\mapsto &\overline{H}\otimes \overline{a}
\end{array}
\end{equation}
 For $u\in\{v^{-1},-v\}$, let $H^{J,u}_w:=H_w\otimes 1 \in \mathcal{L}(u)\otimes_{\He_J}\He$. Denote by $\W^J$ the set of minimal length representatives of $\W/\W_J$. 

\begin{theor}[cf. \cite{Deo87}]\label{Deo87}\leavevmode
\begin{enumerate}
\item For all $w\in \W^J$ there exists a unique element $\underline{H}^{J,v^{-1}}_w\!\!\in\! \mathbf{M}^J$ such that
 \begin{enumerate}
\item $\overline{\underline{H}^{J,v^{-1}}_w}=\underline{H}^{J,v^{-1}}_w$, and 
\item $\underline{H}^{J,v^{-1}}_w=\sum_{y\in \W^J}{m^J_{y,w} H^{J,v^{-1}}_y}$,
\end{enumerate}
where the $m^J_{y,w}$ are such that $m^J_{w,w}=1$ and $m^J_{y,w}\in v \mathbb{Z}[v]$ if $y\neq w$.
\item For all $w\in \W^J$ there exists a unique element $\underline{H}^{J,-v}_w\!\!\in\! \mathbf{N}^J$ such that
\begin{enumerate}
\item $\overline{\underline{H}^{J,-v}_w}=\underline{H}^{J,-v}_w$, and 
\item $\underline{H}^{J,-v}_w=\sum_{y\in \W^J}{n^J_{y,w}H^{J,-v}_y}$,
\end{enumerate}
\end{enumerate}
where the $n^J_{y,w}$ are such that $n^J_{w,w}=1$ and $n^J_{y,w}\in v \mathbb{Z}[v]$ if $y\neq w$.
\end{theor}

  
\begin{rem}\label{rem_parab_regularHM}In the case $J=\0$, the two parabolic modules coincide with the regular module: $\mathbf{M}^{\0}=\mathbf{N}^{\0}=\He$. Moreover, $\underline{H}^{\0,v^{-1}}_w=\underline{H}^{\0,-v}_w=\underline{H}_w$ for all $w\in W$.
\end{rem}

From now on, we will focus on the case $u=v^{-1}$, that is we will deal only with $\mathbf{M}^J$. 
 The action of the Hecke algebra $\He$ on $\MeJ$ is defined as follows. Let $s\in\SR$ be a simple reflection and let $x\in\W^J$, then we have (cf. \cite[\S3]{Soe97}):

\begin{equation}\label{eqn_multPH} \underline{H}_s\cdot H_x^{J,v^{-1}} =\left\{
\begin{array}{lcl}
H_{sx}^{J,v^{-1}}+vH_{x}^{J,v^{-1}}&\text{\small{if}}& sx\in\W^J, sx>x,\\
H_{sx}^{J,v^{-1}}+ v^{-1}H_{x}^{J,v^{-1}}&\text{\small{if}}& sx\in\W^J, sx<x,\\
(v+v^{-1})H_{x}^{J,v^{-1}}&\text{\small{if}}& sx\not\in\W^J.
\end{array}
\right.
\end{equation}

\subsection{Definition of the categorification of $\MeJ$}\label{ssec_categorification} For any category  $\mathscr{C}$ which is exact in the sense of Quillen (cf. \cite{Qu}), let us denote by $[\mathscr{C}]$ its Grothendieck group; that is, the abelian group with generators 
\begin{equation*}[X],\qquad \hbox{\text{for}}\quad X\in \text{Ob}(\mathscr{C}),
\end{equation*} and relations 
\begin{equation*}
[Y]=[X]+[Z]\qquad\hbox{\text{for every exact sequence}}\quad 0\rightarrow X \rightarrow Y \rightarrow Z\rightarrow 0.
\end{equation*} 

For an exact endofunctor $F$ on $\mathscr{C}$, denote by $[F]$ the 
induced endomorphism of $[\mathscr{C}]$.

By a \emph{categorification} of $\MeJ$, we mean an exact category $\mathscr{C}$ together with an autoequivalence $G$ and a 
family of exact endofunctors $\{F_s\}_{s\in\SR}$ satisfying the following requirements:
\begin{itemize}\label{Def_catPHM}
\item[(C1)] $[\mathscr{C}]$ becomes an $\mathcal{L}$-module via $v^{i}\cdot [A]=[G^i A]$ for any $i\in\mathbb{Z}$ and there is
an isomorphism $h^J: [\mathscr{C}]\stackrel \sim \longrightarrow \MeJ$ of $\mathcal{L}$-modules; 
\item[(C2)]  for any simple reflection $s\in\SR$, we have an isomorphism of functors $GF_s\cong F_s G$;
\item[(C3)]  for any simple reflection $s\in \SR$, the following diagram commutes:

\begin{displaymath}
\begin{xymatrix}{
[\mathscr{C}]\ar[d]_{h^J}\ar[r]^{[F_s]}&[\mathscr{C}]\ar[d]^{h^J}\\
\MeJ \ar[r]_{ \underline{H}_s \cdot}&\MeJ.}
\end{xymatrix}
\end{displaymath}

\end{itemize}

\begin{rem} Our notion of  $\MeJ$-categorification differs from the one of Mazorchuk and Stroppel (cf. \cite{MS1}, Remark 7.8). 
Indeed, we made the (weaker) requirement of $\mathscr{C}$ being  exact instead of abelian. If we take the above categorification, restrict it to the additive category of projective objects and then abelianise it in the standard way, then this abelianisation is a 2-functor  (see \cite[\S3.3]{M}) and will transform the above categorification into a categorification using abelian categories, in the spirit of \cite{MS1}. 
\end{rem}

\begin{rem}
In \cite{W}, Williamson studied the 2-category of singular Soergel bimodules. A full tensor subcategory of it (${}^{\0}\mathcal{B}^J$ in his notation) also provides a categorification of $\MeJ$.
\end{rem}

The main goal of this paper is to construct such a categorification. In particular, we will generalise a categorification
 of the Hecke algebra obtained by Fiebig in \cite{Fie07a}, which is known, by results in \cite{Fie08b}, to be equivalent to the one via Soergel's bimodules in \cite{Soe07}.

\section{Sheaves on moment graphs}

\subsection{Moment graphs}

In this section we recall some definitions from \cite{Fie08a}, \cite{Fie08b}, \cite{FieNotes}.

\begin{defin}[cf. {\cite{Fie08b}}]\label{Def_MG} Let  $k$ be a field,  let $V$ be a finite dimensional $k$-vector space and $\mathbb{P}(V)$ the corresponding projective space. A $V$\emph{-moment graph} is given  by a tuple $(\V,\E, \trianglelefteq, l)$ where:
\begin{itemize}
\item[(MG1)] $(\V,\E)$ is a  graph with a set of vertices $\V$ and a set of edges $\E$;
\item[(MG2)] $ \trianglelefteq$ is a partial order on $\V$ such that $x,y\in \V$ are comparable if they are  linked by an edge;
\item[(MG3)] $l:\E\rightarrow \mathbb{P}(V)$  is a map  called the \emph{label function}.
\end{itemize}
\end{defin}

\begin{rem}This is the traditional definition (see \cite{Fie08b}). We note that the fact that $\V$ is equipped with a partial order (similarly to the notion of quasi-hereditary algebra)  is used only in the definition of Braden-MacPherson sheaves.
\end{rem}

As in \cite{Fie08b}, we think of the order as giving each edge a direction: we write $E:x\rightarrow y\in \E$ if $x\leq y$. We write $x-\!\!\!-\!\!\!-y$ or $y-\!\!\!-\!\!\!-x$  if we want to ignore the order.

\subsubsection{Bruhat graphs}

Let $(\W, \SR)$ be a Coxeter system and denote by $m_{st}$ the order of the product of two simple reflections $s,t\in \SR$.
Let $V$ be the geometric representation of $(\W, \SR)$ (cf. \cite[\S5.3]{Humph}). Then  $V$ is a real vector space with basis indexed by the set of simple reflections $\Pi=\{\als \}_{s\in\SR}$ and $s$ acts on $V$  by 
\begin{equation*}s: v\mapsto v-2\langle v,\als\rangle \als.
\end{equation*}
where $\langle\cdot, \cdot\rangle:V\times V\longrightarrow \mathbb{R}$ the symmetric bilinear form given by
\begin{equation*}
\langle\alpha_s,\alpha_t\rangle=
\left\{
\begin{array}{lcl}
-\cos\left(\frac{\pi}{m_{st}}\right)&\text{if}&m_{st}\neq \infty,\\
-1&\text{if}&m_{st}=\infty.
\end{array}
\right.
\end{equation*}

Consider a subset $J\subseteq \SR$ and keep the same notation as in the previous section. Choose $\lambda\in V$ such that $\W_J=\text{Stab}_{\W}(\lambda)$ . Then $\W^J$ can be identified with the orbit $\W\cdot\lambda$ via $x\mapsto x(\lambda)$.

Recall that the set of reflections  $\T$ of $\W$  is
\begin{equation*}\T=\left\{ ws w^{-1}\vert\ s\in\SR,\ w\in\W\right\}.
\end{equation*}

\begin{defin}[cf. {\cite[\S2.2]{Fie08b}}]\label{defin_BruhatMG}The \emph{Bruhat moment graph} $\MGJ$ associated to the Coxeter datum $(\W, \SR, J)$ is the following $V$-moment graph:
\begin{itemize}
\item the set of vertices is given by $\W\cdot \lambda\!\leftrightarrow\!\W^J$ and $x\rightarrow y$ is an edge and only if $\ell(x)<\ell(y)$ and there exists a reflection $t\in \T$ such that  $x(\lambda)=ty(\lambda)$, that is $y=txw$, for some $w\in\W_J$, and $y\not\in x\W_J$.
\item the partial order $W^J$ is   the (induced) Bruhat order;
\item $l(x\rightarrow txw)$ is given  by the line generated by $x(\lambda)-tx(\lambda)$ in $\mathbb{P}(V)$.
\end{itemize}
\end{defin}

Consider now  two Bruhat moment graphs on  $V$: $\MG=\MG(\W, \SR, \0)$ and $\MGJ=\MG(\W, \SR, J)$. The \emph{canonical quotient map}  $p^J:\MG\rightarrow \MGJ$ is induced by the map  $p^J_{\V}:x\rightarrow x^J$, with $x^J$ minimal length representative of the coset $x\W_J$.

\begin{expl}Let $\W=S_3$, the symmetric group on three letters. In this case we have $V=\mathbb{R}^2$, $\Pi=\{\alpha, \beta\}$, and the angle between the two roots is $\frac{2\pi}{3}$.  If we  fix $J=\{s_{\alpha}\}$, then $p^J$ is as follows.

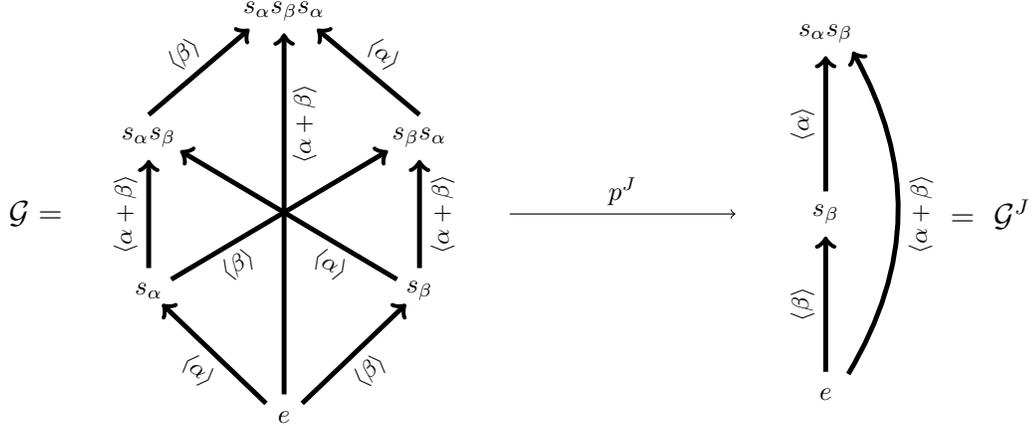
\vspace{3mm}
\begin{tikzpicture}[scale=0.6]
\node at (-0.5,4)  {\Large{$\mathcal{G}=$}};
\node at (5,-0.5)  {$e$};
\node at (8,2.3)  {$s_{\beta}$};
\node at (2,2.3)  {$s_{\alpha}$};
\node at (5,8.5)  {$s_{\alpha}s_{\beta}s_{\alpha}$};
\node at (8,5.7)  {$s_{\beta}s_{\alpha}$};
\node at (2,5.7)  {$s_{\alpha}s_{\beta}$};
\path[->,line width=0.00001mm] (5,4) edge node[below, sloped]{\small{$\langle\alpha+\beta\rangle $}} (5,8);
\path[->,line width=0.65mm] (4.6,-0.2) edge node[below, sloped]{\small{$\langle\alpha\rangle$}} (2.3,2);
\path[->,line width=0.65mm] (5.4,-0.2) edge node[below, sloped]{\small{$\langle\beta\rangle$}} (7.7,2);
\path[->,line width=0.65mm] (2,6.2) edge node[above, sloped]{\small{$\langle\beta\rangle $}} (4.23,8.1) ;
\path[->,line width=0.65mm] (8,6.2) edge node[above, sloped]{\small{$\langle\alpha\rangle $}} (5.73,8.1);
\path[->,line width=0.65mm] (5,0) edge (5,8);
\path[->,line width=0.65mm] (2,2.8) edge node[above, sloped]{\small{$\langle\alpha+\beta\rangle $}} (2,5.2);
\path[->,line width=0.65mm] (8,2.8) edge node[below, sloped]{\small{$\langle\alpha+\beta\rangle$}} (8,5.2);
\path[->,line width=0.65mm] (7.5,2.55) edge (2.7,5.4);
\path[line width=0.00001mm] (7.5,2.55) edge node[below, sloped]{\small{$\langle\alpha\rangle $}} (5,4);
\path[->,line width=0.65mm] (2.5,2.55) edge (7.3,5.4);
\path[-,line width=0.00001mm] (2.5,2.55) edge node[below, sloped]{\small{$\langle\beta\rangle $}} (5,4);
\path[->,line width=0.15mm] (10,4) edge node[above]{{$p^J$}} (15,4);
\node at (17,0)  {$e$};
\node at (17,4)  {$s_{\beta}$};
\node at (17,8)  {$s_{\alpha}s_{\beta}$};
\path[->,line width=0.65mm] (17,0.5) edge node[above,sloped]{\small{$\langle\beta\rangle$}}  (17,3.5);
\path[->,line width=0.65mm] (17,4.5) edge node[above, sloped]{\small{$\langle\alpha\rangle$}}   (17,7.5);
\path[->,line width=0.65mm] (17.5,0.45) edge[bend right] node[below, sloped]{\small{$\langle\alpha+\beta\rangle$}}  (17.6,7.6);
\node at (20.65,4)  {\Large{$=\,\,\mathcal{G}^{J}$}};
\end{tikzpicture}\\ We have $p_{\V}^J(e)=p_{\V}^J(s_{\alpha})=e$, $p^J_{\V}(s_{\beta})=p^J_{\V}(s_{\beta}s_{\alpha})=s_{\beta}$ and $p_{\V}^J(s_{\alpha}s_{\beta})=p_{\V}^J(s_{\alpha}s_{\beta}s_{\alpha})=s_{\alpha}s_{\beta}$. 
\end{expl}

\subsection{Sheaves on a $V$-moment graph}

\subsubsection{Conventions} For any finite dimensional vector space $V$ over the field $k$ (with ${\rm char}\ k\neq 2$), 
 we denote by $S=\text{Sym}(V)$ its symmetric algebra. Then $S$ is a polynomial ring and we provide it with 
the grading induced by setting $S_{\{2\}}=V$. From now on, all the $S$-modules will be finitely generated and $\mathbb{Z}$-graded. Moreover, we will consider only degree zero morphisms between them. For a graded $S $-module $M=\oplus_i M_{\{i\}}$ and for $j\in \mathbb{Z}$, we denote by $M\langle j\rangle$ the $\mathbb{Z}$-graded $S$-module obtained from $M$ by shifting the grading by $j$, that is $(M\langle j\rangle)_{\{i\}}=M_{\{j+i\}}$.

\begin{defin}[cf. \cite{BM01}]\label{Def_ShMG} Let $\MG=(\V, \E, \tle, l)$ be a $V$-moment graph, then a \emph{sheaf } $\F$ \emph{on }$\MG$ is given by 
$(\{\F^x\}, \{\F^E\},$ $\{\rho_{x,E}\})$, where 
\begin{itemize}
\item[(SH1)] for all $x\in \V$, $\F^x$  is an $S$-module;
\item[(SH2)] for all $E\in\E$, $\F^{E}$ is an $S$-module such that $l(E)\cdot \F^E=\{0\}$;
\item[(SH3)] for $x\in \V$, $E\in \E$, $\rho_{x,E}:\F^x\rightarrow \F^E$ is  a homomorphism of $S$-modules defined if $x$ is incident to the edge $E$.
\end{itemize}
\end{defin}

\begin{rem} We may consider the following topology on the space $\Gamma=\V\cup\E$ (cf. \cite[\S1.3]{BM01}).
We say that a subset $O\subseteq \Gamma$ is open, if whenever a vertex $x$ is in $O$, then all edges adjacent to 
$x$ are also in $O$. With this topology, the object in Definition \ref{Def_ShMG}  is a sheaf of $S$-modules on $\Gamma$ in the usual sense.
For our purposes, it will be sufficient to consider sheaves as purely combinatorial, algebraic objects.
\end{rem}

\begin{expl}[cf. {\cite[\S1]{BM01}}]\label{Def_StrSh} Let $\MG=(\V, \E, \tle, l)$  be a $V$-moment graph, then its \emph{structure sheaf} $\ShZ$ is given by
\begin{itemize}
\item for all $x\in\V$, $\ShZ^x=S$;
\item for all $E\in\E$, $\ShZ^E=S/l(E)\cdot S$;
\item    for all $x\in \V$ and $E\in \E$, such that $x$ is incident to the edge $E$, $\rho_{x,E}:S\rightarrow S/l(E)\cdot S$ is the canonical quotient map.
\end{itemize}
\end{expl}

\begin{defin}[cf. \cite{FieNotes}]\label{Def_morphShMG} Let $\MG=(\V, \E,\tle, l)$  be a $V$-moment graph and let $\F=(\{\F^x\}, \{\F^E\}, \{\rho_{x,E}\})$, $\F'=(\{\F'^x\}, \{\F'^E\}, \{\rho'_{x,E}\})$ be two sheaves on $\MG$.  A \emph{morphism} $\varphi:\F\longrightarrow \F'$  is given by the following data:
\begin{itemize}
\item[(MSH1)] for all $x\in\V$, $\varphi^x:\F^{x}\rightarrow {\F'}^x$ is a homomorphism of $S$-modules;
\item[(MSH2)] for all $E\in\E$,  $\varphi^E:\F^{E}\rightarrow {\F'}^E$ is a homomorphism of $S$-modules
such that, if $x\in\V$ is  incident to the edge $E$, the following diagram commutes:
\end{itemize}
\begin{equation*}
\begin{xymatrix}{
\F^x \ar[d]_{\varphi^x}\ar[r]^{\rho_{x,E}}&\F^E \ar[d]^{\varphi^E}\\
{\F'}^x \ar[r]_{\rho'_{x,E}}&{\F'}^E.\\
}
\end{xymatrix}
\end{equation*}
\end{defin}

\begin{defin}\label{Def_catShMG} Let $\MG$  be a $V$-moment graph. We denote by $\ShMGk$ the category of sheaves
 on $\MG$ and corresponding morphisms.
\end{defin}

\begin{rem}
Observe that the category of sheaves on $\MG $ is graded,  with the shift of grading autoequivalence  $\langle 1 \rangle : \ShMGk\rightarrow \ShMGk$ given by 
$$(\{\F^x\}, \{\F^E\}, \{\rho_{x,E}\})\mapsto(\{\F^x\langle 1 \rangle\}, \{\F^E\langle 1 \rangle\}, \{\rho_{x,E}\circ\langle 1 \rangle\}).$$
Moreover $\ShMGk$ is an additive category, with zero object $(\{0\}, \{0\}, \{0\})$ and biproduct given by
\begin{equation*}
(\{\F^x\}, \{\F^E\}, \{\rho_{x,E}\})\oplus (\{\F'^x\}, \{\F'^E\}, \{\rho'_{x,E}\})=(\{\F^x\oplus\F'^x\}, \{\F^E\oplus\F'^E\}, \{(\rho_{x,E},\rho'_{x,E})\}),
\end{equation*}
and idempotent split.
\end{rem}

\subsection{Sections of a sheaf on a moment graph}\label{ssec_secSh}

Even if $\ShMGk$  is not a category of sheaves in the usual sense, we may define the notion of sections following \cite{Fie08a}.

\begin{defin}Let $\MG=(\V,\E, \tle, l)$ be a $V$-moment graph,  $\F=(\{\F^x\}, \{\F^E\}, \{\rho_{x,E}\})\in\ShMGk$ and $\mathcal{I}\subseteq \V$. Then the \emph{set of sections of} $\F$ \emph{over} $\mathcal{I}$ is denoted by $\Gamma(\I,\F)$ and defined as
\begin{equation*}\Gamma(\mathcal{I},\F):=\left\{ (m_x)\in \prod_{x\in \mathcal{I}} \F^{x} \  \Big\vert\  \begin{array}{c} 
\rho_{x,E}(m_x)=\rho_{y,E}(m_y)\\ \text{for all } E:x-\!\!\!-\!\!\!- y\in\E ,\, x,y\in\mathcal{I}                                                                                                                  
                                                                                                                 \end{array}
   \right\}.
\end{equation*}
\end{defin}

We set $\Gamma(\F):=\Gamma(\V,\F)$, that is the set of \emph{global sections} of $\F$. 


\begin{expl} A very important example is given by the set of global sections of the structure sheaf 
$\ShZ$ (cf. Example \ref{Def_StrSh}). In this case, we get the \emph{structure algebra}:
\begin{equation*}\Z:=\Gamma(\ShZ)=\left\{(z_x)_{x\in\V}\in\prod_{x\in\V}S \ \Big\vert\  \begin{array}{c} z_x-z_y\in l(E)\cdot S \\  \text{for all } E:x-\!\!\!-\!\!\!- y\in\E
\end{array}
\right \} .
\end{equation*}
\end{expl}

\begin{rem}The algebra $\Z$ should be thought of as the center of a non-critical block in the deformed category $\mathcal{O}$ (cf. \cite[Theorem 3.6]{Fie}).
\end{rem}

It is easy to check that $\Z$, equipped with componentwise addition and multiplication, is an algebra and that there is an action of $S$ on it by diagonal multiplication. Moreover, for any sheaf $\F\in\ShMGk$, the structure algebra $\Z$ acts on the space $\Gamma(\F)$ via componentwise multiplication, so $\Gamma$ defines a functor from the category of sheaves on $\MG$ to the category of $\Z$-modules:
\begin{equation}\Gamma:\ShMGk\rightarrow \Z\text{-mod}.
\end{equation}

\subsection{BMP-sheaves}\label{sec_BMP_sh} Let $\MG=(\V, \E,\tle, l)$  be a $V$-moment graph. For all $\F\in \ShMGk$ and $x\in \V$, we set
\begin{equation*}\E_{\delta x}:=\left\{  E\in\E\mid \text{there is } y\in\V   \text{ with } E:x\rightarrow y         \right\},
\end{equation*}
 \begin{equation*}\V_{\delta x}:=\left\{ y\in \mathcal{V}\mid \text{there is } E\in\E_{\delta x}  \text{ with } E:x\rightarrow y       \right\}.
   \end{equation*}

Additionally, for any $x\in\V$ denote $\{ \triangleright x \}=\{ y\in \V\,\,|\,\, y\triangleright x   \}$ and define $\F^{\delta x}$ as the image of $\Gamma(\{ \triangleright x \}, \F )$ under the composition of the following functions:
\begin{equation*}
u_x:\Gamma\left(\{\triangleright x\},\F\right)\longrightarrow \bigoplus_{y\triangleright  x}\!\F^y\longrightarrow\bigoplus_{y\in \V_{\delta x}}\!\!\F^y\stackrel{\oplus \rho_{y,E}}{\longrightarrow}\bigoplus_{E\in \E_{\delta x}} \!\!\F^E
.
\end{equation*}

\begin{theor}[cf. \cite{BM01}]\label{BM01}Let $\MG=(\V, \E,\tle, l)$  be a $V$-moment graph and let $w\in\V$. There exists a unique up to isomorphism indecomposable sheaf $\BMP(w)$ on $\MG$ with the following properties:  
\begin{itemize}
\item[(BMP1)] If $x\in \V$, then $\BMP(w)^x\cong 0$, unless $x\trianglelefteq w$. Moreover, $\BMP(w)^w\cong S$;
\item[(BMP2)] If $x,y\in \V$, $E:x\rightarrow y\in \E$, then the map  $\rho_{y,E}:\BMP(w)^y\rightarrow \BMP(w)^E$ is surjective with kernel  $l(E)\cdot\BMP(w)^y$;
\item[(BMP3)]  If $x,y\in \V$, $x\neq w$ and $E:x\rightarrow y\in E$, then $\rho_{\delta x}:=\bigoplus_{E\in \E_{\delta x}}\rho_{x,E}: \BMP(w)^x\rightarrow \BMP(w)^{\delta x}$ is a projective cover in the category of graded $S$-modules.
\end{itemize}
\end{theor}

We call $\BMP(w)$  the \emph{BMP-sheaf}.

\section{Modules over the structure algebra}

Let $\Z$ be the structure algebra (see \S\ref{ssec_secSh}) of a regular Bruhat graph $\MG=\MG(\W,\0)$ and denote by $\Z\text{-mod}^{\text{f}}$ the  category of $\mathbb{Z}$-graded $\Z$-modules that are torsion free and finitely generated  over $S$. In \cite{Fie08b},
 Fiebig defined translation functors on  the category  $\Z\text{-mod}^{\text{f}}$. Using these, he  defined inductively a full subcategory 
$\H$ of $\Z$-mod and proved that $\H$, in characteristic zero, is equivalent to a category of bimodules introduced by Soergel in \cite{Soe07}. 
In \cite{Fie07a} it is shown that $\H$ categorifies the Hecke algebra $\He$ (and the periodic module $\mathbf{M}$), using translation functors. The aim of this chapter is to define translation functors in the parabolic setting and to extend some results of \cite{Fie07a}.  
 
Let $\W$ be a Weyl group, let $\SR$ be its set of simple reflections and let $J\subseteq \SR$. Hereafter we will keep the notation
 we used in \S\ref{sec_HeckeM}. Recall that, for any $z\in\W$, there is a unique factorisation $x=x^Jx_J$, with $x^J\in\W^J$, 
 $x_J\in\W_J$ and $\l(x)=\l(x^J)+\l(x_J)$ (cf. \cite[Proposition 2.4.4]{BB}).

In \cite{Fie08b}, for all $s\in \SR$, an involutive automorphism $\sigma_s$ of the structure algebra of a regular Bruhat graph is
 given. In a similar way, we will define an involution $\ssig$ for a fixed simple reflection $s\in \SR$ on the structure algebra 
$\ZJ$ of the parabolic Bruhat  moment graph $\MGJ$. 

Let $x,y\in\W^J$. Notice that $l(x-\!\!\!-\!\!\!-y)=\al_t$ if and only if $l(sx-\!\!\!-\!\!\!-sy)=s(\al_t)$, because
 $sxw(sy)^{-1}=sxwy^{-1}s=sts$, for some $w\in\W_J$. 

Denote by  $\tau_s$ the automorphism of the symmetric algebra $S$ induced by the mapping $\lambda\mapsto s(\lambda)$ for all $\lambda\in V$. For any $(z_{x})_{x\in\W^J}\in \ZJ$, we set $\ssig\big((z_{x})_{x\in\W^J}\big)=(z_{x}')_{x\in\W^J}$, where $z_{x}':=\tau_s(z_{(sx)^J})$. This is again an element of the structure algebra from what we have observed above. 

Let us fix the following notation:

\begin{itemize}
\item $\sZJ$ denotes  the space of invariants with respect to ${}_s\sigma$;
\item $\msZJ$ denotes the space of anti-invariants with respect to ${}_s\sigma$. 
\end{itemize}

 We denote moreover by $\overline{\als}$ the element of $\ZJ$ whose  components are all equal to $\als$. We obtain the following decomposition of $\ZJ$ as a $\sZJ$-module:

\begin{lem}\label{lem_left_invariants}$\ZJ= \sZJ\oplus \overline{\al_s}\cdot \sZJ$.
\end{lem}
\proof Because  $\ssig$ is  an involution
, we get $\ZJ= \sZJ\oplus \msZJ$. Since $\overline{\als}\in  \ZJ$ and   $s(\als)=-\als$, it follows $\ssig(\overline{ \als})=-\overline{\als}$ and so  $\overline{\als}\cdot \sZJ\subseteq\! \msZJ$ and we now have to prove the other inclusion, that is every element $z\in \!\msZJ$ is divisible  by $\overline{\als}$ in $ \msZJ$.

If $z=(z_{x})\in \msZJ$, then, for all $x\in\W^J$, 
\begin{equation*}z_{x}=-\tau_s(z_{(sx)^J})\equiv -z_{(sx)^J} \,\,\,\,\,(\text{mod } \als).
\end{equation*} 
On the other hand, 
\begin{equation*}
z_{x}\equiv z_{(sx)^J} \,\,\,\,\,(\text{mod } \als).
\end{equation*}
 It follows that $2z_{x}\equiv 0$  $(\text{mod } \als)$, that is $\als$ divides $z_{x}$ in $S$. 

It remains to verify that $z':=\overline{(\als)^{-1}}\cdot z\in \Z$, that is   $z'_{x}-z'_{(tx)^J}\equiv 0 \,(\text{mod } \alt)$ for any $x\in \W^J$ and $t\in\T$. If $(tx)^J=(sx)^J$ there is nothing to prove, since $\als$ divides $z'_{x}=z_{(sx)^J}$ and $z'_{(sx)^J}=z_{x}$, and hence also their difference. On the other hand, if $(tx)^J\neq (sx)^J$ we get the following:
\begin{equation*}\als\cdot(z'_{x}-z'_{(tx)^J})=z_{x}-z_{(tx)^J}\equiv 0\, (\text{mod } \alt).
\end{equation*}
Since $\als$ and $\alt$ are linearly independent, $\als\not \equiv 0\,(\text{mod } \alt)$ and we obtain 
\begin{equation*}z'_{x}-z_{(tx)^J}'\equiv 0\,(\text{mod } \alt). \qedhere
\end{equation*}
\endproof

\subsection{Translation functors and special modules}\label{sssec_TransLeft}

In order to define translation functors, we need an action of $S$ on $\sZJ$ and $\ZJ$. 

\begin{lem}\label{lem_clJ}For any $\lambda\in V$ and any $x\in\W^J$, let us set
\begin{equation}c(\lambda)^J_x:=\sum_{x_J\in\W_J} xx_J(\lambda).
\end{equation}
Then $c(\lambda)^{J}:=(c(\lambda)^J_x)_{x\in\W^J}\in\sZJ$.
\end{lem}
\proof First we prove that $c(\lambda)^{J}\in\ZJ$,
 that is $c(\lambda)^J_x- c(\lambda)^J_{(tx)^J}\equiv 0\,(\text{mod } \alt)$. Since for any $x_J$ there exists an element $y_J$ such that
 $xx_J=t\,(tx)^J\,y_J$, we obtain
\begin{align*}\sum_{x_J\in\W_J} xx_J(\lambda)-\sum_{x_J\in\W_J}(tx)^J\, x_J(\lambda)&=\sum_{y_J\in\W_J} t\,(tx)^J\,y_J(\lambda)-\sum_{y_J\in\W_J} (tx)^J\,y_J(\lambda)\\
&=t\left( \sum_{y_J\in\W_J} \,(tx)^J\,y_J(\lambda)\right)-\sum_{y_J\in\W_J} (tx)^J\,y_J(\lambda)\\
&= \left(\sum_{y_J\in\W_J}2 \left\langle (tx)^J\,y_J(\lambda), \alt\right\rangle\right ) \alt\\
&\equiv 0\,(\text{mod }\al_t)\ .
\end{align*}

To conclude it is left to show that $c(\lambda)^J$ is invariant with respect to ${}_s\sigma$. For any $x\in \W^J$, one has
\begin{align*}
\tau_s\left( c(\lambda)^J_{x}\right)&=\tau_s\left( \sum_{x_J\in\W_J} xx_J(\lambda)\right)\\
&=\sum_{x_J\in\W_J}sxx_J(\lambda)\\
&=c(\lambda)^J_{sx}.
\end{align*}
Hence we have ${}_s\sigma(c(\lambda)^J)=(\tau_s c(\lambda)^J_{sx})_{x\in\W^J}=c(\lambda)^J$.
\endproof

For any $x\in\W^J$, denote by $\eta_{x}$ the endomorphism of the symmetric algebra $S$ induced by the map $\lambda\mapsto c(\lambda)_x^J $ 
for all $\lambda\in V$.  Now, by  Lemma \ref{lem_clJ}, the action of $S$  on $\ZJ$ given by \begin{equation}\label{eqn_actionS}p.(z_{x})_{x\in\W^J}=(\eta_x(p)z_x)\,\,\,\,\,p\in S\,,\,\,z\in\ZJ,
 \end{equation} preserves $\sZJ$. Thus any $\ZJ$-module or $\sZJ$-module has an $S$-module structure as well.
 Let  $\ZJ\text{-mod}^{\text{f}}$, resp.  $\sZJ\text{-mod}^{\text{f}}$, be the category of $\mathbb{Z}$-graded  $\ZJ$-modules, resp. $\sZJ$-modules,  that are torsion free and finitely generated over $S$.

The \emph{translation on the wall} is the functor ${}^{s,on}\theta:\ZJ$ -mod $\rightarrow\sZJ$-mod defined by the mapping $M\mapsto \text{Res}_{\ZJ}^{\sZJ} M$.

The \emph{translation out of the wall}  is the functor ${}^{s,out}\theta:\sZJ$  -mod $\rightarrow\ZJ$-mod defined by the mapping 
$N\mapsto \text{Ind}_{\ZJ}^{\sZJ} N=\ZJ\otimes_{\sZJ}N$. Observe that this functor is well-defined due to Lemma \ref{lem_left_invariants}.

By composition, we get a functor $\sthJ:={}^{s,out}\theta\circ{}^{s,on}\theta: \ZJ$-mod$\rightarrow \ZJ$-mod that we call \emph{(left) translation functor}.

\begin{rem}This construction is very similar to the one in \cite{Soe90}, where translation functors are defined in the finite case for the coinvariant algebra.
\end{rem}

\begin{rem}
One could consider the idempotent split additive tensor category generated by the translation functors we defined above and describe indecomposable projective. This would be useful in order to strengthen our main result to a proper categorification (see Remark \ref{rem_2cat}). In this paper we are not going to investigate this category of translation functors, but the one of special modules, defined in \S\ref{ssec_HJDef}.
\end{rem}

The following proposition describes the first properties of $\sth$:

\begin{prop}\leavevmode
\begin{enumerate}
 \item The functors from $\sZJ$-mod to $\ZJ$-mod mapping $M\mapsto \ZJ\langle2\rangle\otimes_{\sZJ} M$ and $M\mapsto Hom_{\sZJ}(\ZJ, M)$ are isomorphic.
 \item  The functor $\sth=\ZJ\otimes_{\sZJ} -:\ZJ\text{-mod}\rightarrow\ZJ\text{-mod}$ is selfadjoint up to a shift.
\end{enumerate}
\end{prop}
\proof 
(1) Let $M\in \sZJ$-mod,  we want to prove that $\ZJ\langle2\rangle\otimes_{\sZJ} M\cong Hom_{\sZJ}(\ZJ, M)$ as $\ZJ$-modules. 

First, we show that $\ZJ\langle2\rangle\cong \Hom_{\sZJ}(\ZJ,\sZJ)$ as $\sZJ$-modules. By Lemma \ref{lem_left_invariants}, $\{\overline{1}, \overline{\alpha_s}\}$ is a $\sZJ$-basis for $\ZJ$. Let $\overline{1}^*, \overline{\alpha_s}^*\in\Hom_{\sZJ}(\ZJ,\sZJ)$ be the $\sZJ$-basis  dual to $\overline{1}$ and $\overline{\alpha_s}$, that is 
\begin{equation*}
\overline{1}^*(\overline{1})=\overline{1},\quad \overline{1}^*(\overline{\alpha_s})=\overline{0}, \quad\overline{\alpha_s}^*(\overline{\alpha_s})=\overline{1},\quad \overline{\alpha_s}^*(\overline{1})=\overline{0},
\end{equation*}
where $\overline{1}\in\sZJ$, resp. $\overline{0}\in\sZJ$, is thesection with 1, resp. 0, in all entries. Since $\deg(1)-2=-2=\deg(\overline{\alpha_s}^*)$ and $\deg(\overline{\alpha_s})-2=0=\deg{\overline{1}^*}$, we have an isomorphism of $\sZJ$-modules $\ZJ\langle2\rangle\cong \Hom_{\sZJ}(\ZJ,\sZJ)$ defined by the mapping 
\begin{equation*}\overline{1}\mapsto \overline{\alpha_s}^*, \quad \overline{\alpha_s}\mapsto \overline{1}^*.
\end{equation*}
 
Because $\ZJ$ is free of rank two over $\sZJ$,  $\Hom_{\sZJ}(\ZJ, M)\cong\Hom_{\sZJ}(\ZJ, \sZJ)\otimes_{\sZJ}M $ by the map
\begin{equation*}
\varphi\mapsto \overline{\alpha_s}^*\otimes \varphi(\overline{\alpha_s})+\overline{1}^*\otimes \varphi(\overline{1})	 .
\end{equation*} This conclude the proof of (1).

(2) Since $\ZJ\otimes_{\sZJ}-$ and $\Hom_{\sZJ}(\ZJ,-)$ are, resp.,  left and right adjoint to the restriction functor, we obtain the following chain of isomorphisms for any pair $M,N\in\ZJ$:
\begin{align*}
\Hom_{\ZJ}(\sth M, N)&= \Hom_{\ZJ}\left(\ZJ\otimes_{\sZJ}(\text{Res}_{\ZJ}^{\sZJ} M), N\right)\\
&\cong \Hom_{\ZJ}\left(\text{Res}_{\ZJ}^{\sZJ} M, \text{Res}_{\ZJ}^{\sZJ} N\right)\\
&\cong \Hom_{\ZJ}\left(M, \text{Hom}_{\sZJ}(\ZJ,\text{Res}_{\ZJ}^{\sZJ} N)\right)\\
&\cong \Hom_{\ZJ}\left(M,\ZJ\langle 2\rangle \otimes_{\sZJ}(\text{Res}_{\ZJ}^{\sZJ} N)\right)\\
&=\Hom_{\ZJ}( M,\sth \langle 2 \rangle N).\qedhere
\end{align*}
\endproof



\subsection{Parabolic special modules}\label{ssec_HJDef}

As in \cite{Fie08b}, we define, inductively, a full  subcategory of $\ZJ$-mod.

Let $B^J_e\in\ZJ$-mod be the free $S$-module of rank one on which $z=(z_{x})_{x\in \W^J}$ acts via multiplication by $z_e$.

\begin{defin} \leavevmode
 \begin{itemize}
\item The category  $\HJ$ of \emph{special $\ZJ$-modules} is the full subcategory of $\ZJ\text{-mod}^{\text{f}}$ whose objects  are isomorphic to a direct summand of a direct sum of modules of the form ${}^{s_{i_1}}\theta\circ\ldots \circ{}^{s_{i_r}}\theta(B^J_e)\langle n\rangle$, where $s_{i_1}, \ldots, s_{i_r}\in\SR$ and $n\in\mathbb{Z}$.
 \item The category $\sHJ$ of \emph{special $\sZJ$-modules} is the full subcategory  of $\sZJ\text{-mod}^{\text{f}}$ whose objects are isomorphic to a direct summand of ${}^{s,on}\theta(M)$ for some $M\in \HJ$.
\end{itemize}
\end{defin}

\subsection{Finiteness of special modules}

Let $\Omega$ be a finite subset of $\W^J$ and denote by $\ZJ(\Omega)$ the sections of the structure sheaf over $\Omega$, that is 

\begin{equation*}\ZJ(\Omega)=\left \{(z_{x})\in\prod_{x\in \Omega} S\,\Big\vert 
\begin{array}{c}
z_{x}\equiv z_{y} \,\,\,(\!\!\!\!\mod \alt)\\
\text{ if there is } w\in\W_J \text{ s.t. }y\,w\,x^{-1}=t\in\T 
\end{array}\right \}.
\end{equation*}

If $\Omega \subseteq \W^J$ is $s$-invariant,  that is $s\Omega=\Omega$, we may restrict $\ssig$ to it. We denote by $\sZJ(\Omega)\subseteq \ZJ(\Omega)$ the space of invariants and, using Lemma \ref{lem_left_invariants}, we get a decomposition $\ZJ(\Omega)=\sZJ(\Omega)\oplus \overline{\als} \cdot \sZJ(\Omega)$.

In the following lemma we prove, the \emph{finiteness} of special $\ZJ$-modules, as Fiebig does in \cite{Fie07a} 
for special $\Z$-modules.

\begin{lem}\label{lem_finPSM}\leavevmode
\begin{enumerate}
\item Let $M\in \HJ$. Then there exists a finite subset $\Omega\subset \W^J$ such that $\ZJ$ acts on $M$ via the canonical map $\ZJ\rightarrow \ZJ(\Omega)$.
\item Let $s\in \SR$ and let $N$ be an object in $\sHJ$. Then there exists a finite $s$-invariant subset $\Omega\subset \W^J$  such that $\sZJ$ acts on $N$ via the canonical map $\sZJ\rightarrow \sZJ(\Omega)$.
\end{enumerate}
\end{lem}
\proof We prove (1) by induction. 
It holds clearly for $B_e$, since $\ZJ$ acts on it via the map $\ZJ\rightarrow\ZJ(\{e\})$. Now we have to show 
 that if the claim is true for $M\in\HJ$, then it holds also for $\sth(M)$. 
Suppose $\ZJ$ acts via the map $\ZJ\rightarrow \ZJ(\Omega)$ over $M$. Observe that we may assume $\Omega$ $s$-invariant, 
since we can just replace it by $\Omega\cup s\Omega$, which is still finite. In this way the $\sZJ$-action on $\sth M$ factors via $\sZJ\rightarrow\sZJ(\Omega)$ and so we obtain $\sth M:=\ZJ\otimes_{\sZJ}M=\ZJ(\Omega)\otimes_{\sZJ(\Omega)}M$.

Claim (2) follows directly from claim (1).
\endproof

\section{Modules with Verma flag and statement of the main result}

We recall some notation from \cite{Fie08a}. Let $Q$ be the quotient field of $S$ and let $A$ be an  $S$-module. Then we
 denote by $A_{Q}=A\otimes_{S}Q$. Let us assume $\MG$ to be such that  for any $M\in\Z-\text{mod}^{\text{f}}$
there is a canonical decomposition 
\begin{equation}\label{decomp_MQ}
M_{Q}=\bigoplus_{x\in \V} M_{Q}^{x}
\end{equation} and so a canonical inclusion $M\subseteq \bigoplus_{x\in \V} M_{Q}^{x}$. For all subsets  of the set of vertices $\Omega\subseteq \V$, we may define:
\begin{equation*}M_{\Omega}:=M\cap \bigoplus_{x\in \Omega} M_{Q}^{x},\end{equation*}
\begin{equation*}M^{\Omega}:=M/M_{\V\setminus \Omega}=\im\left( M\rightarrow M_{Q}\rightarrow\bigoplus_{x\in \Omega} M_{Q}^{x}\right).
\end{equation*}

 For all $x\in\V$, we set
\begin{equation*}M_{[x]}:=\ker\left(M^{ \{\trianglerighteq x\}}\rightarrow M^{ \{\triangleright x\}}\right )
\end{equation*}
and, if $x\triangleleft y$ and $[x,y]=\{x,y\}$, we denote 
\begin{equation*}M_{[x,y]}:=\ker\left(M^{ \{\trianglerighteq x\}}\rightarrow M^{ \{\trianglerighteq x\}\setminus\{x, y\}}\right).
\end{equation*}

\begin{rem} In \cite{Fie08a} the module $M_{[x]}$ is denoted by $M^{[x]}$. The notation we are adopting in this paper is the one of \cite{Fie07a}.
\end{rem}

\subsection{Modules with a Verma flag}
 
From now on, let $\MG$ be a Bruhat moment graph. In \cite{Fie08a} it is shown that in this case any $M\in \Z\text{-mod}^{\text{f}}$ admits a decomposition like (\ref{decomp_MQ}) and hence the modules $M_{[x]}$ are well defined for any $x\in\V$.

Let $\mathscr{V}$ denote the full subcategory of $\Z\text{-mod}^{\text{f}}$ whose objects admit a 
Verma flag, that is $M\in\mathscr{V}$ if and only if $M^{\Omega}$ is a graded free $S$-module for any $\Omega\subseteq \V$ upwardly closed with respect to the partial order in the set of vertices. In our hypotheses this condition is equivalent to $M_{[x]}$ being a graded free $S$-module for any $x\in \V$ (cf. \cite[Lemma 4.7]{Fie08a}).

\subsubsection{Exact structure}\label{sssec_exact_str}
Now we want to equip the category $\mathscr{V}$ with an exact structure .

\begin{defin}Let $L\rightarrow M\rightarrow N$ be a sequence in  $\mathscr{V}$. We say that it is \emph{short exact} if 
\begin{equation*}0\rightarrow L_{[x]}\rightarrow M_{[x]}\rightarrow N_{[x]} \rightarrow 0
\end{equation*}
is a short exact sequence of $S$-modules for any $x\in\V$.
\end{defin}
\begin{rem}This is not the original definition of exact structure Fiebig gave in \cite{Fie08a}, which was on the whole category $\Z\text{-mod}^{\text{f}}$, but it is
 known to be equivalent to it if we only consider the category $\mathscr{V}$, 
that is precisely the one we are dealing with (cf. \cite[Lemma 2.12]{Fie08b}).
\end{rem}

\subsection{Decomposition and subquotients of modules on $\bf{\ZJ}$}\label{ssec_sec_subquotPSM}

Lemma \ref{lem_sthMx} describes the action of $\sth$ on the subquotients $M_{[x]}$, for $x\in \V$. This is important in order to show that
 $\HJ$ categorifies the parabolic Hecke algebra. Before proving Lemma \ref{lem_sthMx}, 
we need a combinatorial result, that follows easily from the so--called lifting lemma, which we now recall.

\begin{lem}[``Lifting lemma'', cf. {\cite[ Lemma 7.4]{Humph}}]\label{lift_lem} Let $s\in \SR$ and $v,u\in \W$ be such that $vs<v$ and $u<v$.
\begin{enumerate}
\item If $us<u$, then $us<vs$;
\item if $us>u$, then $us\leq v$ and $u\leq vs$.
\end{enumerate}
Thus, in both cases, $us\leq v$.
\end{lem}

\begin{lem}\label{lem_combpar1}Let $x\in\W^J$ and $t\in \SR$. If $tx\not\in\W^J$, then $(tx)^J=x$.
\end{lem}
\proof If $tx\not \in \W^J$, then there exists a simple reflection $r\in J$ such that $txr<tx$ and, since $x\in\W^J$, $xr>x$. Using (the left version of) Lemma \ref{lift_lem} (1)  with $s=t$, $v=xr$ and $u=tx$, we get $txr<x$. Applying Lemma \ref{lift_lem} (1) with $s=r$, $v=x$ and $u=txr$ it follows $tx>x$. Finally, from Lemma \ref{lift_lem} (2) we obtain $txr\leq x$, that, together with $x<xr$, gives $txr=x$.
\endproof

\begin{lem}\label{lem_sthMx}Let $s\in\SR$ and $x\in\W^J$, then
\begin{equation*}({}^s\theta M)_{[x]}\cong \left\{
\begin{array}{lcl}
M_{[x]}\langle-2\rangle\oplus M_{[sx]}\langle-2\rangle&\text{\small{if}}& sx\in\W^J, sx> x,\\
M_{[x]}\oplus M_{[sx]}&\text{\small{if}}& sx\in\W^J, sx<x,\\
M_{[x]}\langle-2\rangle\oplus M_{[x]}&\text{\small{if}}& sx\not\in\W^J.
\end{array}
\right.
\end{equation*}
\end{lem}
\proof By Lemma \ref{lem_combpar1}, if $sx\not\in\W^J$, then $(sx)^J=x$ and $M_{[x]}\in\sZJ$-mod, so by Lemma 
\ref{lem_left_invariants} we get $\ZJ\otimes_{\sZJ}M_{[x]}=M_{[x]}\langle-2\rangle\oplus M_{[x]}$.

If $x\neq sx$, we have a short exact sequence  of $S$-modules $0\rightarrow M_{[x]}\rightarrow M_{[x,sx]} \rightarrow M_{[sx]}\rightarrow 0$. By Lemma \ref{lem_left_invariants} the module $\ZJ$ is flat over $\sZJ$, which is a finitely generated free S-module. Hence we have  $\sth M_{[x,sx]}=\ZJ\otimes_{\sZJ} M_{x,sx}=(\sth M)_{[x,sx]}=\sth M_{[x]}\oplus \sth M_{[sx]}$.
 Moreover $\sth M_{[x,sx]}=\ZJ(\{x,sx\})\otimes_{\sZJ(\{x,sx\})} M_{[x,sx]}$ and the two isomorphisms follow keeping in mind that
 $\ZJ(\{x,sx\})_{[x]}\cong S\langle-2\rangle$ if $x< sx$, while $\ZJ(\{x,sx\})_{[x]}\cong S$ if $x> sx$. 
\endproof

Using induction, we obtain  the following corollary:

\begin{cor}\label{cor_freeSM}Let $M\in\HJ$, then for any $x\in\W^J$, $M_{[x]}$ is a finitely generated torsion free $S$-module ad hence $M\in \mathscr{V}$.
\end{cor}

In this way we get an exact structure also on $\HJ$ and we are finally able to state the main result of this paper:

\begin{theor}\label{Thm_catHM} The category $\H^J$ together with the shift in degree $\langle -1 \rangle$ and (shifted) translation functors is
 a categorification of the parabolic Hecke module $\MeJ$.
\end{theor}

\begin{rem}\label{rem_2cat} Theorem \ref{Thm_catHM} could be strengthen to a proper categorification by presenting the result as a 2-representation of a 2-category. The 2-category to be considered is the one generated by the translation functors we defined in \S\ref{sssec_TransLeft} and the 2-representation to look at is given by the action  of these functors on the category $\HJ$ we constructed in \S\ref{ssec_HJDef}. The question of describing indecomposable 1-morphisms in this category, which we are not going to address in this paper, seems seems to be very interesting. 
\end{rem}

\begin{rem} The recent results of Elias and Williamson (cf. \cite{EW}) imply that the results of \cite{MS1} transfer to all Coxeter systems.
\end{rem}

\section{Proof of the categorification theorem} The proof of Theorem \ref{Thm_catHM} consists of several steps:
\begin{enumerate}
\item we show that the functor $\sth\circ\langle 1 \rangle$ is exact (Lemma \ref{sth_exact});
\item we define the character map $h^J:[\HJ]\rightarrow \MeJ$ (\S\ref{ssec_charmaps});
\item we observe that the map $[\langle  -1 \rangle]: [\HJ]\rightarrow [\HJ]$ provide $[\HJ]$ with a structure of $\mathcal{L}$-module and that $h^J$ is a map of $\mathcal{L}$-modules  (\S\ref{ssec_charmaps});
\item via explicit calculations, we prove that the functors $\sth\circ \langle 1\rangle$, $s\in \SR$,  satisfy (C3), that is the maps they induce on $[\HJ]$ commute with $h^J$  (Proposition \ref{prop_left_mult});
\item we demonstrate that the character map is surjective by choosing a certain basis for $\MeJ$ and showing that every element of this basis has a preimage in $[\HJ]$ under $h^J$ (Lemma \ref{Lem_surj_hJ});
\item  we prove that the character map is surjective (Lemma \ref{Lem_inj_hJ}) using a description of indecomposable special modules in terms of Braden-MacPherson sheaves (Proposition \ref{prop_HJ_BMP}).
\end{enumerate}
This concludes the proof, since (C2), that is $\langle  -1 \rangle\circ( \sth\circ \langle 1\rangle) \cong (\sth\circ \langle 1\rangle)\circ  \langle  -1 \rangle$ for any $s\in \SR$, is trivially satisfied.

We start by proving the exactness of  shifted translation functors.

\begin{lem}\label{sth_exact}For any $s\in \SR$ the functor ${}^s\theta\circ\langle 1\rangle:\HJ\rightarrow\HJ$ is exact. 
\end{lem}
\proof
Let $L\rightarrow M\rightarrow N$ be an exact sequence, then for any $x\in\V$ 
\begin{equation*}0\rightarrow L_{[x]}\rightarrow M_{[x]}\rightarrow M_{[x]} \rightarrow 0
\end{equation*}
is a short exact sequence of $S$-modules. In particular, also
\begin{equation*}0\rightarrow L_{[sx]}\rightarrow M_{[sx]}\rightarrow N_{[sx]} \rightarrow 0
\end{equation*}
is short exact.
The claim follows immediately from Lemma \ref{lem_sthMx}, and the fact that finite direct sums and shifts preserve exactness.
\endproof

\subsection{Character maps}\label{ssec_charmaps}
Let $A$ be a $\mathbb{Z}$-graded,  free and finitely generated $S$-module; then  $A\cong\bigoplus_{i=1}^n S\langle k_i\rangle$, for some $k_i\in\mathbb{Z}$. We can associate to $A$ its \emph{graded rank}, that is the following Laurent polynomial:
\begin{equation*}\rk A:=\sum_{i=1}^n v^{-k_i}\in \mathcal{L}.
\end{equation*} 
This is well-defined, because the $k_i$'s are uniquely determined, up to reordering.

Let $M\in \HJ$, then by Corollary \ref{cor_freeSM}, we may define a map $h^J:[ \HJ ]\rightarrow \MeJ$ as follows.
\begin{equation*}h^J([M]):=\sum_{x\in\W^J}v^{\l(x)}\rk M_{[x]} \ H_x^{J,v^{-1}}\in\MeJ.
\end{equation*}

The Grothendieck group $[\HJ]$ is equipped with a structure of $\mathcal{L}$-module via $v^i[M]=[M\langle -i\rangle]$. 
Observe that  for any $M\in \HJ$ one has $h^{J}(v[M])=h^{J}([M\langle -1 \rangle])=vh^J([M])$ and so  $h^J$ is a map of $ \mathcal{L}$-modules.

\begin{prop}\label{prop_left_mult} For  each $M\in\HJ$ and for any $s\in\mathcal{S}$ we have $h^J([\sth M\langle1\rangle])=\underline{H}_s\cdot h^J([M])$, that is the following diagram is commutative:
\begin{displaymath}
\begin{xymatrix}{
[\HJ]\ar[d]_{h^J}\ar[r]^{[\sth\circ \langle 1 \rangle]}&[\HJ]\ar[d]^{h^J}\\
\MeJ \ar[r]_{ \underline{H}_s \cdot}&\MeJ.}
\end{xymatrix}
\end{displaymath}
\end{prop} 
\proof By Lemma \ref{lem_sthMx}, for any $x\in\W^J$ we have 
\begin{equation*}\rk({}^s\theta M)_{[x]}= \left\{
\begin{array}{lcl}
v^2\left(\rk M_{[x]}+\rk M_{[sx]}\right)&\text{\small{if}}& sx\in\W^J, sx>x,\\
\rk M_{[x]}+\rk M_{[sx]}&\text{\small{if}}& sx\in\W^J, sx<x,\\
(v^2+1) \rk M_{[x]}&\text{\small{if}}& sx\not\in\W^J.
\end{array}
\right.
\end{equation*}
Then,
\begin{align*}
h^J([\sth M\langle1\rangle])&=\sum_{x\in\W^J}v^{\l(x)-1}\rk (\sth M)_{[x]} H_x^{J,v^{-1}}\\
&=\sum_{\substack{x\in\W^J, sx\in\W^J\\sx>x}}v^{\l(x)+1}\left(\rk M_{[x]}+  \rk M_{[sx]}\right)H_x^{J,v^{-1}}\\
&\quad +\sum_{\substack{x\in\W^J, sx\in\W^J\\sx<x}}v^{\l(x)-1}\left(\rk M_{[x]}+  \rk M_{[sx]}\right)H_x^{J,v^{-1}}\\
&\quad +\sum_{\substack{x\in\W^J, sx\not\in\W^J}}(v^{\l(x)+1}+v^{\l(x)-1})\rk M_{[x]}H_x^{J,v^{-1}}.
\end{align*}
Finally,
\begin{align*}
\underline{H}_s\cdot h^J([M])&=\sum_{x\in\W^J}v^{\l(x)}(\rk M_{[x]}) \underline{H}_s\cdot H_x^{J,v^{-1}} \\
&=\sum_{\substack{x\in\W^J sx\in\W^J\\sx>x}}v^{\l(x)}(\rk M_{[x]})  (H_{sx}^{J,v^{-1}}+v H_x^{J,v^{-1}})\\
&\quad +\sum_{\substack{x\in\W^J, sx\in\W^J\\sx<x}}v^{\l(x)}(\rk M_{[x]}) (H_{sx}^{J,v^{-1}}+v^{-1}H_x^{J,v^{-1}})\\
&\quad +\sum_{\substack{x\in\W^J, sx\not\in\W^J}}v^{\l(x)}\rk M_{[x]}(v+v^{-1})H_x^{J,v^{-1}}\\
&=\sum_{\substack{x\in\W^J, sx\in\W^J\\sx>x}}\Big[(v^{\l(x)} v \,\rk M_{[x]})+(v^{\l(sx)}\rk M_{[sx]})\Big]H_x^{J,v^{-1}}\\
&\quad +\sum_{\substack{x\in\W^J, sx\in\W^J\\sx<x}}\Big[(v^{\l(x)} v^{-1} \,\rk M_{[x]})+(v^{\l(sx)}\rk M_{[sx]})\Big]H_x^{J,v^{-1}}\\
&\quad +\sum_{\substack{x\in\W^J, sx\not\in\W^J}}(v^{\l(x)+1}+v^{\l(x)-1})\rk M_{[x]}H_x^{J,v^{-1}}\\
&=h^J([\sth M\langle 1\rangle]).
 \qedhere\end{align*}
\endproof

\subsection{The character map is an isomorphism}

In order to prove that $(\H^J, \langle -1 \rangle, \{\sth\circ \langle 1\rangle\})$ is a categorification of $\MeJ$, 
the only step left is to show that $h^J$ is an isomorphism.

\begin{lem}\label{Lem_surj_hJ} The map $h^J:[\H^J]\rightarrow \MeJ$ is surjective.
\end{lem}
\proof We start by defining a  basis of $\MeJ$. Let us set $\widetilde{\underline{H}}^{J,v^{-1}}_e=\underline{H}^{J,v^{-1}}_e$. For any $x\in\W^J$ with $\l(x)=r>0$, let a fix a reduced $x=s_{i_1}\ldots s_{i_r}$, with $s_{i_1}, \ldots, s_{i_r}\in \SR$, denote
\begin{equation*}
\widetilde{\underline{H}}^{J,v^{-1}}_x=\underline{H}^{J,v^{-1}}_{s_1}\cdot\ldots\cdot\underline{H}^{J,v^{-1}}_{s_r}.
\end{equation*}
From Theorem \ref{Deo87}, it follows 
\begin{equation}
\widetilde{\underline{H}}^{J,v^{-1}}_x=H^{J,v^{-1}}_{x}+\sum_{\substack{y\in\W^J\\y< x}} p_{y}H^{J,v^{-1}}_{y}, \qquad 	\hbox{\text{with}} \quad p_z\in\mathbb{Z}[v,v^{-1}].
\end{equation}

Since the set $\{H^{J,v^{-1}}_x\}_{x\in \W^J}$ is a basis of $\MeJ$ as a $\mathbb{Z}[v,v^{-1}]$-module, also $\{\widetilde{\underline{H}}^{J,v^{-1}}_x\}_{x\in \W^J}$ is a basis for $\MeJ$ and it is enough to show that, for any $x\in \W^J$, there exists an object $H\in \H^J$ such that $h^J([H])=\widetilde{\underline{H}}^{J,v^{-1}}_x$.

By definition, $h^J(B^J_e)=M_e=\underline{H}^{J,v^{-1}}_e$. By applying Proposition \ref{prop_left_mult}, we obtain 
\begin{equation*}
 h^J({}^{s_{i_1}}\theta\circ \ldots\circ  {}^{s_{i_r}}\theta B^J_e\langle n\rangle)= (\underline{H}^{J,v^{-1}}_{s_1}\cdot\ldots\cdot\underline{H}^{J,v^{-1}}_{s_r})M_e=
\underline{H}^{J,v^{-1}}_{s_1}\cdot\ldots\cdot\underline{H}^{J,v^{-1}}_{s_r}=\widetilde{\underline{H}}^{J,v^{-1}}_{x}.
\end{equation*}
This conclude the proof of the lemma. \qedhere
\endproof

Proposition \ref{prop_HJ_BMP} will allow us to see any element in $\HJ$ as the space of global sections of some BMP-sheaf on $\MGJ$. From now on, we will denote by $B^J(w)$ the space of global sections of the indecomposable BMP-sheaf $\BMP^J(w)\in\Sh{\MGJ}$. Let us recall a fundamental characterisation of $B^J(w)$.

\begin{theor}[cf. {\cite[Theorem 5.2.]{Fie08b}}] \label{thm_BMP_proj}For any $w\in\MG^J$, the module  $B^J(w)\in \mathscr{V}$ is indecomposable 
and projective. Moreover, every indecomposable projective object in $\mathscr{V}$   is isomorphic to $B^J(w)\langle k \rangle$ for a unique $w\in\MG^J$ and a unique $k\in\mathbb{Z}$.
\end{theor}

\begin{prop}\label{prop_HJ_BMP}A module $M\in \ZJ\text{-mod}^{\text{f}}$ is an indecomposable special module  if and  only there exist a BMP-sheaf $\BMP\in\Sh{\MGJ}$ and $k\in\mathbb{Z}$ such that $M\cong \Gamma(\BMP\langle k \rangle)$ as $\ZJ$-modules.
\end{prop}
\proof By induction, from the exactness of $\sthJ$, it follows that the objects of   $\HJ$ are all projective and then, by Theorem \ref{thm_BMP_proj}, any $M\in\HJ$ may be identified (up to a shift) with the space of global sections of a BMP-sheaf on $\MGJ$. 

We want now to show  that, for any $x\in\W_J$, $B^J(x)\in\HJ$. We prove the claim by induction on $\#\  \text{supp}(M)$, where $\text{supp}(M)=\{x\in\W^J\,\vert\, M^x\neq 0\}$. Clearly, $B_e\cong B^J(e)$.

The statement follows straightforwardly, once proved that, if $sx>x$,  then $\sthJ( B^J(x))=B^J(sx)\oplus B$.

First we show that $\text{supp}(\sthJ ( B^J(x)))\subseteq \{\leq sx\}$, that is $(\sthJ B^J(x))^y=0$ for all $y\not \in   \{\leq sx\}\cap \W^J$. From  Lemma \ref{lem_sthMx}, it follows easily  that $(\sthJ (B^J(x)))_{[y]}=0$ for all $y\not \in   \{\leq sx\}\cap \W^J$. 

Let us observe that, as $\sthJ B^J(x)\in\HJ$, from what we have proved above, there exist $w_1,\ldots, w_r\in\W^J$ and $k_1, \ldots, k_r$ such that $\sthJ (B^J(x))=\oplus_{i=1}^r B^J(w_i)\langle k_i\rangle$ and, for any $y\in\W^J$, 
\begin{equation*} \Big(\bigoplus_{i=1}^r B^J(w_i)\langle k_i\rangle\Big)_{[y]}=\bigoplus_{i=1}^r B^J(w_i)_{[y]}\langle k_i\rangle.
\end{equation*}
So, in particular, for all $y\not \in   \{\leq sx\}\cap \W^J$, 
\begin{align*}
0&=B^J(w_i)_{[y]}\\
&=\ker (\rho_{\delta y}:\BMP^J(w_i)^{y}\rightarrow \BMP^J(w_i)^{\delta y}).
\end{align*}
This implies $\BMP^J(w_i)^{y}=B^J(w_i)^y=0$ for all $i=1, \ldots r$, and so
\begin{equation*}\sthJ (B^J(x))=\bigoplus_{i=1}^r B^J(w_i)\langle k_i\rangle 
\end{equation*}
where $w_i\in\{\leq sx\}$ for al $i=1, \ldots , r$.

It is left to show that there exists at least one $i\in\{1, \ldots, r\}$ such that $w_i=sx$. By applying once again Lemma  \ref{lem_sthMx}, we obtain $(\sthJ(B^J(x)))^{sx}=(\sthJ(B^J(x)))_{[sx]}\cong S$ and hence the statement.\endproof

\begin{lem}\label{Lem_inj_hJ}The map $h^J:[\H^J]\rightarrow \MeJ$ is injective.
\end{lem}
\proof
By Theorem \ref{thm_BMP_proj} and Proposition \ref{prop_HJ_BMP} we know that $\left\{[B^J(w)]\right\}_{w\in\W^J}$ is a $\mathbb{Z}[v,v^{-1}]$-basis of $[\HJ]$ and so every element $Y\in [\HJ]$ can be written as $Y=\sum a_w [B^J(w)]$, with $a_x\in\mathbb{Z}[v,v^{-1}]$ . Let us suppose $Y\in \ker(h^J)$, then
\begin{equation*}
0=h^J(Y)=\sum_{w\in\W^J}a_w \sum_{x\in\W^J}v^{\l(x)}\rk B^J(w)_{[x]} \ H_x^{J,v^{-1}}=
\sum_{x\in\W^J}\left(\sum_{w\in\W^J}v^{\l(x)}a_w\  \rk B^J(w)_{[x]}\right) \ H_x^{J,v^{-1}}.
\end{equation*}
Since the elements $\ H_x^{J,v^{-1}}$ are linearly independent, it follows that $\sum_{w\in\W^J}v^{\l(x)}a_w\  \rk B^J(w)_{[x]}=0$ for any $x\in \W^J$.

If  it were $Y\neq 0$, then we would find a maximal element  $\overline{w}$ such that $a_{\overline{w}}\neq 0$. By (BMP1), we obtain $B^J(w)_{[\overline{w}]}=0$ for all $w<\overline{w}$ and  $B^J(\overline{w})_{[\overline{w}]}\cong S$.  Then,
\begin{equation*}
0=\sum_{w\in\W^J}v^{\l(x)}a_w\  \rk B^J(w)_{[\overline{w}]}=v^{\l(x)}a_{\overline{w}}\  \rk B^J(\overline{w})_{[\overline{w}]}=v^{\l(x)}a_{\overline{w}}\rk S=v^{\l(x)}a_{\overline{w}}.
\end{equation*}
The chain of equalities above gives us a contradiction, since we assumed $a_{\overline{w}}\neq 0$.\endproof

This concludes the proof of Theorem \ref{Thm_catHM}.

\section{The functor $I$}
In this section we define an exact functor $I:\H^J\rightarrow \H^{\0}$ such that the following diagram commutes:
\begin{displaymath}\label{diag_I}
 \begin{xymatrix}{
[\HJ]\ar[d]_{h^J}\ar@{^{(}->}[r]^{[I]}&[\H^{\0}]\ar[d]^{h^{\0}}\\
\MeJ \ar@{^{(}->}[r]_{ i }&\He ,}
\end{xymatrix}
\end{displaymath}

where $i:\MeJ\hookrightarrow \He$ is the map of $\mathcal{L}$-modules given by 
\begin{equation}H_x^{J,v^{-1}}\mapsto \sum_{z\in \W_J}v^{\l(w_J)-\l(z)} H_{xz},
\end{equation}
with $w_J$ the longest element of $\W_J$.

The map $i$ is interesting since it gives us a way to see the parabolic Hecke module $\MeJ$ as submodule of $\He$ and hence its categorification tells us that we can think about the $\HJ$ as a subcategory of $\H^{\0}$. 

We construct the functor $I$  by using a localisation-globalisation procedure. More precisely, we first map the elements of $\HJ$ to certain sheaves on $\MGJ$, then apply a pullback functor mapping them to  sheaves on $\MG$ and finally we take global sections of the latter. A priori it is not clear that we obtain an object in $\H^{\0}$. This fact is shown in Lemma \ref{lem_I_HJtoH}. We then demonstrate the exactness of $I$ (Proposition \ref{prop_exactnessI}) and the commutativity of diagram (\ref{diag_I}) (Proposition \ref{prop_commutativityDiagr}) by a study of the subquotients involved in the definition of the character map.  The realisation of special modules in terms of Braden-MacPherson sheaves given  in the previous section (Proposition \ref{prop_HJ_BMP}) plays a crucial role in the  proof of any of the above results.

\subsection{Construction of the functor $I$}
 The definition of $I$ involves Fiebig's localisation functor  $\mathscr{L}$(cf. \cite[\S3.3]{Fie08a}), which allows us to see 
objects of $\ZJ$-mod as sheaves on the parabolic Bruhat moment graph $\MGJ$.

 Let us assume $\MG$ to be such that  for any $M\in\Z-\text{mod}^{\text{f}}$
there is a canonical decomposition like the one in (\ref{decomp_MQ}). Let $\Z$ be the corresponding structure algebra and $M\in\Z\text{-mod}^{\text{f}}$. 
For any vertex $x\in\V$, we set
\begin{equation}\mathscr{L}(M)^x=M^x.
\end{equation}
For any edge $E:x-\!\!\!-\!\!\!-y$, let us consider  $\Z(E)=\{(z_x,z_x)\in S\oplus S\,\vert\, z_x-z_y\in l(E)S\}$ and $M(E):=\Z(E)\cdot M^{x,y}$.  For $m=(m_x,m_y)\in M(E)$, let us set $\pi_x((m))=m_x$, $\pi_y((m))=m_y$. Then we get $\mathscr{L}(M)^E$ as the push--out  in the following diagram of $S$-modules:
\begin{displaymath}
\begin{xymatrix}
{M(E) \ar[r]^{\pi_x}\ar[d]_{\pi_y}&M^x \ar[d]_{\rho_{x,E}}\\
M^y\ar[r]^{\rho_{y,E}}&\mathscr{L}(M)^E . 
}
\end{xymatrix}
\end{displaymath}

This provides us also with the restriction maps $\rho_{x,E}$ and $\rho_{y,E}$.

It is not hard to verify  (cf. \cite[\S3.3]{Fie08a}) that this is a well-defined functor 
\begin{equation}\mathscr{L}:\Z\text{-mod}^{\text{f}}\rightarrow \ShMGk.\end{equation}
 Moreover, 
the localisation functor $\mathscr{L}$ turns out to be left adjoint to $\Gamma$ (cf. \cite[Theorem 3.5]{Fie08a}). Let  
\begin{itemize}
\item $\Z\text{-mod}^{\text{loc}}$  be the full subcategory of $\Z\text{-mod}^{\text{f}}$, whose objects are the elements $M$ such that there is an isomorphism $\Gamma\circ\mathscr{L}(M)\cong M$, and
\item $\ShMGk^{\text{glob}}$ be the full subcategory of $\ShMGk$, whose objects are the elements $\F$ such that there is an isomorphism $\mathscr{L}\circ \Gamma(\F)\cong \F$.
\end{itemize}

\begin{rem}\label{L_preservesBMP} In general, for a given a sheaf $\F$, one has $(\mathscr{L}\circ \Gamma(\F))^x=\Gamma(\F)^x\subseteq \F^x$. If we consider a BMP-sheaf $\BMP$, then by property (BMP3), $\Gamma(\BMP)^x= \BMP^x$ for any vertex $x\in \V$ and  $\mathscr{L}(\Gamma(\BMP))^E\cong \BMP^E$ for any edge $E\in\E$. Therefore $\mathscr{L}\circ \Gamma(\BMP)\cong \BMP$ and $\BMP\in \ShMGk^{\text{glob}}$.
\end{rem}
Thus, the functors $\mathscr{L}$ and $\Gamma$ induce two inverse equivalences:
\begin{displaymath} 
\begin{xymatrix}{
\Z\text{-mod}^{\text{loc}}\ar[r]&\ShMGk^{\text{glob}}\ar[l]
}.
\end{xymatrix}
\end{displaymath}

Let us focus again on the Bruhat case and consider the functor $p^{J,*}: \textbf{Sh}(\MG^J)\rightarrow  \textbf{Sh}(\MG^{\empty})$ defined as follows:

\begin{itemize}
\item for all $x\in\W$, $(p^{J,*}\F)^x:=\F^{x^J}$;
\item for all $E:x-\!\!\!-\!\!\!-y\in\E$\begin{equation*}(p^{J,*}\F)^{E}=\left\{ \begin{array}{ll}
\F^{f_{\V}(x)}/l(E) \F^{f_{\V}(x)} & \text{if } x^J=y^J,\\
\F^{f_{\E}(E)} &\text{otherwise };
\end{array}
\right.
\end{equation*}
\item for all $x\in \W$ and $E\in \E$, such that $E:x-\!\!\!-\!\!\!-y$,\begin{equation*}(p^{J,*}\rho)_{x,E}=\left\{ \begin{array}{ll}
\text{canonical quotient map} & \text{if } x^J=y^J,\\
\rho_{f_{\V}(x),f_{\E}(E)} &\text{otherwise }.
\end{array}
\right. 
\end{equation*}
\end{itemize}

Finally, we set $I:=\langle  - \l(w_J)\rangle \circ \Gamma\circ p^{J,*}\circ \mathcal{L}$. 

 In order to prove that the functor $I$ maps $\HJ$ to $\H$, we need to recall the moment graph analogue of a theorem by Deodhar relating parabolic Kazhdan-Lusztig polynomials and regular ones. 
The following is a reformulation of  Theorem 6.1 of \cite{L11}:

\begin{theor}\label{thm_DeodMG} Let $J\subseteq \SR$ be such that $W_J$ is finite, with longest element $w_J$. Let $w\in \W^J$, then $p^{J,*}(\BMP^J(w))\cong\BMP^{\0}(ww_J) $ as sheaves on $\MG=\MG(\W, \0)$.
\end{theor}

\begin{lem}\label{lem_I_HJtoH}The functor $I$ maps $\HJ$ to $\H$.
\end{lem}
\proof Let $M\in\HJ$, then, by Proposition \ref{prop_HJ_BMP}, there exist $w_1, \ldots w_r\in\W^J$ and $m_1, \ldots m_r\in\mathbb{Z}$
such that $M=\bigoplus_{i=1}^r B^J(w_i)\langle m_i\rangle $. Then, we get the following:
\begin{align*}
I(M)&=I\left(\bigoplus_{i=1}^r B^J(w_i)\langle m_i\rangle\right)\\
&=\bigoplus_{i=1}^r  \Gamma\circ p^{J,*}\circ\mathscr{L}(B^J(w_i))\langle m_i-\l(w_J)\rangle.
\end{align*}
By Remark \ref{L_preservesBMP}, $\mathscr{L}(B^J(w_i))\cong \BMP^J(w_i)$ for any $i$ and, by Theorem  \ref{thm_DeodMG}, we conclude that 
\begin{equation*}I(M)\cong\bigoplus_{i=1}^r B^{\0}(w_iw_J)\langle m_i-\l(w_J)\rangle.
\qedhere
\end{equation*}
\endproof

\subsubsection{Exactness of $I$}

\begin{lem}\label{lem_PulbackSubModules} 
Let $w\in\W^J$. 
Then,  for all $x\in\W$, 
 $$(\Gamma\circ p^{J,*} \BMP^{J}(w))_{[x]}=\left(\prod_{\substack{y\in\V^{\delta x},\\y\in x\W_J}} \al_y \right) B^J(w)_{[x^J]}, $$
where $\alpha_y$ denotes the label of $x\rightarrow y$.
\end{lem}
\proof For $z\in \W^J$ and $E$ an edge of $\MGJ=\MG(\W,J)$, let us denote by $\rho_{z,E}$ the corresponding restriction map. 
Then, we have the following:


\begin{align*}
(\Gamma \circ p^{J,*} \BMP^{J}(w))_{[x]}&=\bigcap_{y\in\V^{\delta x}} \ker \left( (p^{*,J}\rho)_{x,x\rightarrow y}\right)\\
&= \left(\bigcap_{\substack{y\in\V^{\delta x}\\
y\not\in x\W_J}} \ker (\rho_{x^J,x^J\rightarrow y^J})\right)\,\cap\, \left(\bigcap_{\substack{y\in\V^{\delta x},\\
y\in x\W_J}}\ker \pi_{x,x\rightarrow y}
\right),
\end{align*}
where $\pi_{x,x\rightarrow y}:\BMP^J(w)^{x^J}\rightarrow \BMP^J(w)^{x^J}/\al_y \BMP^J(w)^{x^J}$ is the canonical quotient map 
and $\al_y$ is a generator of $l(x\rightarrow y)$.

Let us observe that, by definition, 

\begin{equation*}
 \bigcap_{\substack{y\in\V^{\delta x}\\ y\not\in x\W_J}} \ker (\rho_{x^J,x^J\rightarrow y^J})= B^{J}(w)_{[x^J]}.
\end{equation*}
 
Moreover, since there is at most one edge adjacent to $x$ labeled by a multiple of $\alpha_y$, the labels of  such edges are pairwise linearly independent and we get
\begin{equation*}
 \bigcap_{\substack{y\in\V^{\delta x},\\y\in x\W_J}} \ker \pi_{x,x\rightarrow y}= \prod_{\substack{y\in\V^{\delta x},\\y\in x\W_J}} \al_y\cdot \BMP^J(w)^{x^J}
\end{equation*}

It follows that
\begin{equation*}
(\Gamma\circ p^{J,*} \BMP^{J}(w))_{[x]}= \left(\prod_{\substack{y\in\V^{\delta x},\\y\in x\W_J}} \al_y \right)\BMP^J(w)_{[x^J]} 
\end{equation*}

This concludes the proof of the lemma.
\endproof

\begin{prop}\label{prop_exactnessI}The functor $I$ is exact, with respect to the exact structure in \S\ref{sssec_exact_str}. 
\end{prop}
\proof  Let us take $M,N\in\HJ$, with $M=\bigoplus_{k}B^J(w_k)\langle m_k \rangle $ and $N=\bigoplus_{l}B^J(w_l)\langle n_l \rangle $.

Let us consider the map $f:L\rightarrow M$ and the induced maps $f_{[x^J]}:M_{[x^J]}\rightarrow N_{[x^J]}$ for any $x^J\in\W^J$. Thanks to Lemma \ref{lem_PulbackSubModules}, it is easy to describe $I(f)_{[x]}$. Namely, if
\begin{equation*}\prod_{\substack{y\in\V^{\delta x}\\y\in x\W_J}} \al_y =\al_{i_1}\cdot\ldots\cdot \al_{i_r},
\end{equation*} we obtain
\begin{equation*}\begin{array}{ccccc}I(f)&:&I(M)_{[x]}&\longrightarrow &I(N)_{[x]}\\
&&(\al_{i_1}\cdot\ldots \cdot\al_{i_r}) m &\mapsto&(\al_{i_1}\cdot\ldots \cdot\al_{i_r} )f_{[x]}(m) 
\end{array}
\end{equation*}

It is clear that, if $0\rightarrow L_{[x]}\rightarrow M_{[x]} \rightarrow N_{[x]}\rightarrow 0$ is a short exact sequence of $S$-modules, then $0\rightarrow (IL)_{[x]}\rightarrow (IM)_{[x]} \rightarrow (IN)_{[x]}\rightarrow 0$ is also exact. \endproof

\subsubsection{Commutativity of the diagram}
The last step missing is the commutativity of Diagram \ref{diag_I}.  Before proving it, we need the following preliminary lemma.

\begin{lem}\label{lem_B0BJ}Let $w\in \W^J$  and let $w_J$ be the longest element of $\W_J$. There is an isomorphism $B^{\0}(ww_J)_{[x]}\cong B^J(w)_{[x^J]}\langle 2\l(x)-2\l(x^J)-2\l(w_J) \rangle$ of graded $S$-modules.
\end{lem}
\proof
 By Theorem \ref{thm_DeodMG}, $\BMP^{\0}(ww_J)\cong p^{J,*} \BMP^{J}(w_J)$ as sheaves on $\MG=\MG(\W, \0)$. It follows that,
 for any $x\in\W$, $B^{\0}(ww_J)_{[x]}\cong (\Gamma\circ p^{J,*} \BMP^{J}(w_J))_{[x]}$ as graded $S$-modules and then, by Lemma \ref{lem_PulbackSubModules},  we obtain

\begin{align*}
B^{\0}(ww_J)_{[x]}&\cong  \left(\prod_{\substack{y\in\V^{\delta x},\\ y\in x\W_J}} \al_y \right)B^J(w)_{[x^J]}\\
&\cong  B^J(w)_{[x^J]}\left\langle 2\cdot \#\  \{y\in\V^{\delta x},y\in x\W_J\} \right\rangle.
\end{align*}

Let $x'=(x^J)^{-1}x\in W_J$. Now, if $\T_J$ is the set of reflections of $\W_J$, 
\begin{align*}
 \#\  \{y\in\V^{\delta x},y\in x\W_J\}&=\#\  \{z\in\W_J\,|\text{ there exists } t\in\T_J \text{ s.t. }z=x't\text{ and }x'<z\}\\
&= \l(w_J)-\l(x')\\
&= \l(w_J)-\l(x)+\l(x^J).\qedhere
\end{align*}
\endproof

Finally, we are able to prove the following proposition, which enable us to embed $\HJ$ in $\H$.

\begin{prop}\label{prop_commutativityDiagr}
 The following diagram is commutative:
\begin{displaymath}\label{diag_I}
 \begin{xymatrix}{
[\HJ]\ar[d]_{h^J}\ar@{^{(}->}[r]^{[I]}&[\H^{\0}]\ar[d]^{h^{\0}}\\
\MeJ \ar@{^{(}->}[r]_{ i }&\He .}
\end{xymatrix}
\end{displaymath}
\end{prop}
\proof 
As $I(\bigoplus_{i\in I} B^J(w_i))=\bigoplus I(B^J(w_i))$, it is enough to prove the statement for the module $B^J(w)$. In this case, we have:

\begin{align*}
 I\left(B^J(w)\right)&=\langle  - \l(w_J)\rangle \circ \Gamma\circ p^{J,*}\circ \mathcal{L} \left(B^J(w)\right)\\
&=\langle  -\l(w_J)\rangle \circ\Gamma\circ p^{J,*}\left(\mathscr{B}^J(w)\right)\\
&\cong \langle - \l(w_J)\rangle \circ\Gamma\left(\mathscr{B}^{\0}(ww_J)\right)\\
&=B(ww_J)\langle  -\l(w_J)\rangle.\\
\end{align*}

Thus, if $B^J(w)_{[x^J]}=\bigoplus_{i\in I_{x^J}}S\langle k_i\rangle$, we get
\begin{align*}
h^{\0}\circ [I]([B^J(w)])&=h^{\0}\left(B^{\0}(ww_J)\langle  \l(w_J)\rangle\right)\\
&=\sum_{x\in \W}v^{-\l(w_J)+\l(x)}\rk B^{\0}(ww_J)_{[x]} H_x\\
\text{ \small{(by Lemma \ref{lem_B0BJ})}}&=\sum_{x\in \W}v^{\l(w_J)+\l(x)}\rk \left(B^{J}(w)_{[x^J]}\langle 2\l(x_J)-2\l(w_J) \rangle\right) H_x\\
&=\sum_{x\in \W}v^{-\l(w_J)+\l(x)}(\sum_{i\in I_{x^J}} v^{-2\l(x_J)+2\l(w_J)-k_i}) H_x\\
&=\sum_{x\in \W}v^{\l(w_J)+\l(x)}(\sum_{i\in I_{x^J}} v^{-2\l(x_J)-k_i}) H_x,
\end{align*}	
where $H_x=H_x^{\0,v^{-1}}$. On the other hand, we have the following:

\begin{align*}
i\circ h^{J}([B^J(w)])&= i\left(\sum_{x^J\in\W^J} v^{\l(x^J)}\rk B^J(w)_{[x^J]} H_{x^J}^{J,v^{-1}}\right)\\
&=\sum_{x^J\in\W^J}\left[v^{\l(x^J)}\left(\sum_{i\in I_{x^J}} v^{-k_i}\right) i(H_{x^J}^{J,v^{-1}})\right]\\
&=\sum_{x^J\in\W^J}\left[v^{\l(x^J)}\left(\sum_{i\in I_{x^J}} v^{-k_i}\right) \left(\sum_{x_J\in\W_J} v^{\l(w_J)-l(x_J)} H_{x^Jx_J} \right)\right]\\
&=\sum_{x^J\in\W^J}\sum_{x_J\in\W_J}\left(\sum_{i\in I_{x^J}} v^{\l(x^J)-k_i+\l(w_J)-\l(x_J)} \right) H_{x^Jx_J}\\
&=\sum_{x\in \W}v^{\l(w_J)+\l(x)}\left(\sum_{i\in I_{x^J}} v^{-2\l(x_J)-k_i}\right) H_x.\qedhere
\end{align*}	
\endproof

\section{Connection with equivariant category $\mathcal{O}$}

In this section we briefly discuss the connection of our results with non-critical blocks in an equivariant version of category $\mathcal{O}$. Our main references are \cite{Fie} and \cite{Fie08a}. 

Let $\mathfrak{g}$ be a complex symmetrisable Kac-Moody algebra and $\mathfrak{b}\supseteq\mathfrak{h}$ its Borel and Cartan subalgebras. The Weyl group $\W$ of $\mathfrak{g}$ naturally acts on $\mathfrak{h}^{\star}$, and  we can consider equivalence classes $\Lambda\in\mathfrak{h}^{\star}/\sim$. An element $\lambda\in \mathfrak{h}^{\star}$ is non-critical if $2(\lambda+\rho, \beta)\not \in \mathbb{Z}(\beta,\beta)$ for any imaginary root $\beta$ and an orbit $\Lambda$ is non-critical if  any $\lambda\in\Lambda$ is non-critical. 

Let us fix a non-critical orbit $\Lambda$ and a weight $\lambda_0\in \Lambda$.  As in Definition \ref{defin_BruhatMG}, we can look at the $\W$-orbit of $\lambda_0$, which gives us a Bruhat moment graph on $\mathfrak{h}^{\star}$. We want  to discuss the representation theoretic content of $\HJ$, where $J$ is in this case given by the set of simple reflections generating $\text{Stab}_{\W}\lambda_0$. Let us denote by  $\MG(\Lambda)$ such a graph.

Let  $S= S(\mathfrak{h})$ be  the symmetric algebra of $\mathfrak{h}$,   $R=S_{(\mathfrak{h})}$ be  its localisation at $0\in \mathfrak{h}^{\star}$, and  $\tau:S\rightarrow R$ be the canonical map. For any $\mu\in\mathfrak{h}^{\star}$ and any $(\mathfrak{g}$-$R)$-bimodule $M$, we define its \emph{$\mu$-weight space} as 
\begin{equation*} 
M_{\mu}=\left\{m\in M\mid H.m=(\lambda(H)+\tau(H))m\ \text{ for any } H\in\mathfrak{h}\right\}.
\end{equation*}
If $\mathfrak{g}$-mod-$R$ denotes the category of $(\mathfrak{g}$-$R)$-bimodules, then the equivariant version of category $\mathcal{O}$ we want to study is the following:
\begin{equation*}
\mathcal{O}_R=\left\{M\in\mathfrak{g}\text{-mod-}R\ \Big\vert\ 
\begin{array}{l}
M \text{ is locally finite as } (\mathfrak{b}\text{-}R)\text{-bimodule},\\
M=\bigoplus_{\mu\in\mathfrak{h}^{\star}}M_{\mu}\\
\end{array}
\right\}.
\end{equation*}

For any $\mu\in\mathfrak{h}^{\star}$ let us consider the ($\mathfrak{h}$-$R$)-bimodule $R_{\mu}$ free of rank one over $R$
 on which  $\mathfrak{h}$ acts via the character $\mu+\tau$. The projection $\mathfrak{b}\rightarrow \mathfrak{h}$ allows us to  consider $R_{\mu}$ as a ($\mathfrak{b}$-$R$)-bimodule and we can now induce to obtain the \emph{equivariant Verma module} of weight $\mu$: $M_{R}(\mu)=U(\mathfrak{g})\otimes_{U(\mathfrak{b})}R_{\mu}$, where  $U(\mathfrak{g})$ and $U(\mathfrak{b})$ are  the  enveloping algebras of $\mathfrak{g}$ and $\mathfrak{b}$, respectively.

Let $\mathcal{M}_{R}$ be the full subcategory of $\mathcal{O}_R$ whose objects admit a finite filtration with subquotients isomorphic to equivariant Verma modules. Since $\mathcal{O}_R$ is abelian and $\mathcal{M}_{R}$ is closed under extensions in $\mathcal{O}_{R}$, the category $\mathcal{M}_{R}$ inherits an exact structure.

For an equivalence class $\Lambda\in\mathfrak{h}^{\star}/\sim$, let $\mathcal{O}_{R,\Lambda}$, resp.  $\mathcal{M}_{R,\Lambda}$, be the full subcategory of $\mathcal{O}_R$, resp. $\mathcal{M}_R$,  consisting of all objects $M$  such that the highest weight of every simple subquotient of $M$ lies in $\Lambda$. Then there are block decompositions, according to the following two results.

\begin{prop}[cf. {\cite[Proposition 2.8]{Fie}}]
The functor 
\begin{equation*}
\begin{array}{rcl}
\prod_{\Lambda\in\mathfrak{h}^{\star}/\sim}\mathcal{O}_{R,\Lambda}&\rightarrow &\mathcal{O}_{R}\\
\{M_{\Lambda}\}&\mapsto&\bigoplus_{\Lambda\in \mathfrak{h}^{\star}/\sim}M_{\Lambda}
\end{array}
\end{equation*}is an equivalence of categories.
\end{prop}

\begin{theor}[cf. {\cite[Theorem 6.1]{Fie08a}}]
The functor 
\begin{equation*}
\begin{array}{rcl}
\prod_{\Lambda\in\mathfrak{h}^{\star}/\sim}\mathcal{M}_{R,\Lambda}&\rightarrow &\mathcal{M}_{R}\\
\{M_{\Lambda}\}&\mapsto&\bigoplus_{\Lambda\in \mathfrak{h}^{\star}/\sim}M_{\Lambda}
\end{array}
\end{equation*}is an equivalence of categories.
\end{theor}

Now it is important to notice that we could have substituted $S$ by the local algebra $R$ in the constructions and definitions we have considered and all the results of this paper would have still worked. Let us denote by $\Z_R$ the $R$-version of the structure algebra of $\MG(\Lambda)$ and by $\mathscr{V}_{R,\Lambda}$ the category of  $\Z_R$-modules admitting a Verma flag. The main result of \cite{Fie08b} is the following one:

\begin{theor}[cf. {\cite[Theorem 7.1]{Fie08a}}] There is an equivalence of exact categories 
\begin{equation*}\mathbb{V}:\mathcal{M}_{R,\Lambda}\rightarrow  \mathscr{V}_{R,\Lambda}.
\end{equation*}
\end{theor}

\subsection{Projective objects} For $\nu\in\Lambda$, let  $\Lambda^{\leq \nu}:=\{\lambda\in\Lambda\mid \lambda\leq \nu\}$. We want to consider a truncated version of $\mathcal{M}_{R, \Lambda}$:
\begin{equation*}
\mathcal{M}_{R, \Lambda^{\leq \nu}}=\left\{M\in\mathcal{M}_{R,\Lambda}\ \Big\vert\ \left(M:M_{R}(\mu)\right)\neq 0\text{  only if } \mu\in\Lambda^{\leq \nu}\right\}.
\end{equation*} 
As a reference for the truncated category $\mathcal{O}$, we address the reader to \cite{RW}, where it was introduced.
 
Denote by $\mathscr{V}_{R,\Lambda^{\leq \nu}}$ the category of sheaves on the moment graph $\MG(\Lambda)^{\leq \nu}$, obtained by restricting the set of vertices of  $\MG(\Lambda)$  to $\Lambda^{\leq \nu}$. By \cite[Proposition 3.11]{Fie06}, the functor $\mathbb{V}$ restricts to a functor $\mathbb{V}^{\leq \nu}:\mathcal{M}_{R,\Lambda^{\leq \nu}}\rightarrow \mathscr{V}_{R,\Lambda^{\leq \nu}}$, which is also an equivalence of categories.


Let $\HJ_R$ denote the $R$-version of the category of special modules, and let $\mathcal{H}^{J}_{R, \Lambda^{\leq \nu}}$ be the subcategory of $\HJ_{R}$ consisting of modules having support on $\MG(\Lambda)^{\leq \nu}$. From Theorem \ref{thm_BMP_proj}, a module  $M\in \mathscr{V}_{R, \Lambda^{\leq \nu}}$ is  indecomposable and projective if and only if there exist a $w\in\Lambda^{\leq \nu}$ and a $k\in\mathbb{Z}$ such that $M\cong B^J(w)\langle k \rangle$ and, by Proposition \ref{prop_HJ_BMP}, there exists one and only one  indecomposable $M\in \mathcal{H}^{J}_{R, \Lambda^{\leq \nu}}$ isomorphic to $B^J(w)$. In summary,

\begin{prop}\label{prop_eqvtRepnTh}
Let $P\in \mathcal{M}_{R, \Lambda^{\leq \nu}}$. Then $P$ is  indecomposable, projective if and only if  $\mathbb{V}P$ is an indecomposable special module.
\end{prop}


 For $\lambda_0$ regular, that is $\text{Stab}_{\W}\lambda_0=\{e\}$, this fact has been already proven by Fiebig (cf. \cite{Fie08b}) and  used in the paper \cite{Fie07a}, where the interchanging between local and global description of the image of the projective modules under $\mathbb{V}$  played a fundamental role.


\section*{Acknowledgements} I would like to acknowledge Peter Fiebig for useful discussions and Winston Fairbairn for  helpful conversations.  I owe many thanks to Michael Ehrig and Ben Salisbury for their careful proofreading.


\begin{thebibliography} {9}



\bibitem[BeBe]{BeBe}A. Beilinson, J. Bernstein, Localisation de $\g$-modules, \textit{C. R. Acad. Sci. Paris S\'er. I Math.} \textbf{292} (1981), no. 1, 15--18.


\bibitem[BGG]{BGG}I. N. Bernstein, I. M. Gelfand, S. I. Gelfand, A certain category of g-modules. \textit{Funkcional. Anal. i Prilozen.} \textbf{10} (1976), no. 2, 1--8.


\bibitem[BB]{BB}A. Bj\"orner, F. Brenti, \textit{Combinatorics of Coxeter Groups}, Graduate Texts in Mathematics, Vol. 231, Springer, New York, 2005.


\bibitem[BMP]{BM01}T. Braden, R. MacPherson, From moment graphs to intersection cohomology, \textit{Math. Ann.} \textbf{321} (2001), no.3, 533--551.


\bibitem[BK]{BK}J-L. Brylinski, M. Kashiwara, D\'emonstration de la conjecture de Kazhdan et Lusztig sur les modules de Verma,
 \textit{C. R. Acad. Sci. Paris S\'er.  A-B} \textbf{27} (1980), no. 6, 373--376.

\bibitem[Deo87]{Deo87}V. Deodhar, On some geometric  aspects of Bruhat orderings II. The parabolic analogue of Kazhdan-Lusztig polynomials, \textit{J. Algebra} \textbf{111} (1987), no. 2, 483--506.

\bibitem[EW]{EW}B. Elias, G. Williamson, The Hodge theory of Soergel bimodules, preprint 2012, arXiv:1212.079.

\bibitem[Fie03]{Fie}P.Fiebig, Centers and translation functors for the category $\mathcal{O}$ over Kac-Moody algebras, \emph{Math. Z.} \textbf{243} (2003), no.4, 689--717.


\bibitem[Fie06]{Fie06}P. Fiebig, The combinatorics of category $\mathcal{O}$ over symmetrizable Kac-Moody algebras, \textit{Transform. Groups} \textbf{11} (2006), no. 1, 29--49. 

\bibitem[Fie08b]{Fie08b}P. Fiebig, The combinatorics of Coxeter categories, \textit{Trans. Amer. Math. Soc.} \textbf{360} (2008), no. 8, 4211--4233.


\bibitem[Fie08a]{Fie08a}P. Fiebig, Sheaves on moment graphs and a localization of Verma flags, \textit{Adv. Math.} \textbf{217} (2008), no. 2,  683--712.

\bibitem[FieNotes]{FieNotes} P. Fiebig, Moment graphs in representation theory and geometry, Skript der  Vorlesung zu Bruhatgraphen im Wintersemester 2008/09, available at http://www.algeo.math.uni-erlangen.de/fileadmin/algeo/users/fiebig/Skripten/Skript\_MomGra.pdf.

\bibitem[Fie11]{Fie07a}P. Fiebig, Sheaves on affine Schubert varieties, modular representations and Lusztig's conjecture, \textit{J. Amer. Math. Soc.} \textbf{24} (2011), no.1, 133--181.


\bibitem[FW]{FW}P. Fiebig, G. Williamson, Parity sheaves, moment graphs and the p-smooth 
locus of Schubert varieties,  to appear in \emph{Ann. Inst. Four.}.


\bibitem[GKM]{GKM}M. Goresky, R. Kottwitz, R. MacPherson, Equivariant cohomology, Koszul duality, and the localization theorem, \textit{Invent. Math.} \textbf{131}, no.1 (1998), 23--83.


\bibitem[Humph]{Humph}J. E. Humphreys, \textit{Reflection Groups and Coxeter Groups}, Cambridge Studies in Advanced Mathematics, 29, Cambridge University Press, Cambridge, 1990.


\bibitem[KL]{KL}D. Kazhdan, G. Lusztig, Representations of Coxeter groups and Hecke algebras, \textit{Invent.  Math.} \textbf{53} (1979), no. 2, 165--184.


\bibitem[KL80]{KL80}D. Kazhdan, G. Lusztig, Schubert varieties and Poincar\'e duality, Geometry of the Laplace Operator. ( Proc. Sympos. Pure Math. Univ. Hawaii Honolulu, Hawaii 1979), 185--203, Proc. Sympos. Pure Math. XXXVI  \textit{Amer. Math. Soc.} (1980).



\bibitem[KT]{KT} M. Kashiwara, T. Tanisaki, Kazhdan-Lusztig conjecture for affine Lie algebras with negative level, \textit{Duke Math. J.} \textbf{77} (1995), no. 1, 21--62.


\bibitem[L11]{L11}M. Lanini, Kazhdan-Lusztig combinatorics in the moment graph setting, \textit{J. Algebra} \textbf{370} (2012), 152--170.



\bibitem[M]{M}V. Mazorchuk, Lectures on algebraic categorification, QGM Matsre Class Series, Europen Mathematical Society (EMS), Z\"{u}rich, 2012.

\bibitem[MS1]{MS1}V. Mazorchuk, C. Stroppel, Translation and shuffling of projectively presentable modules and a categorification  of a parabolic Hecke module, \textit{Trans. Amer. Math. Soc.} \textbf{357} (2005), no. 7, 2939--2973.

\bibitem[MS2]{MS2}V. Mazorchuk, C. Stroppel, Categorification of (induced) cell modules and the rough structure ofgeneralized Verma modules, \textit{Adv. Math.} \textbf{219} (2008), no. 4, 1363--1426.

\bibitem[Q]{Qu}D. Quillen, Higher algebraic $K$-theory. I,  Algebraic $K$-theory. I: Higher K-theories, (Proc. Conf. Battelle Memorial Inst., Seattle, Wash., 1972), pp.   85--147, Lecture Notes in Math. Vol. 341, Springer, Berlin 1973.

\bibitem[RW]{RW}A. Rocha-Caridi, R. Wallach, Projective Modules over Graded Lie Algebras, I, \textit{Math. Z.} \textbf{180} (1982), no. 2, 151--177.

\bibitem[Soe90]{Soe90}W. Soergel,  Kategorie $\mathcal{O}$, perverse Garben und Moduls \"{u}ber den Koinvarianten zur Weylgruppe, \textit{J. Amer. Math. Soc.} \textbf{3} (1990), 421--445. 

\bibitem[Soe97]{Soe97}W. Soergel, Kazhdan-Lusztig polynomials and combinatoric[s] for tilting modules, \textit{Represent. Theory} \textbf{1} (1997), 83--114.


\bibitem[Soe07]{Soe07}W. Soergel, Kazhdan-Lusztig Polynome und unzerlegbare Bimoduln \"uber Polynomringen, \textit{J.  Inst.  Math.  Jussieu} \textbf{6} (2007), no. 3, 501--525. 

\bibitem[W]{W}G. Williamson, Singular Soergel Bimodules, \textit{Int. Math. Res. Not.}, IMNR 2011, no. 20, 4555--4632.

\end{thebibliography}
\end{document}